\documentclass[a4paper,12pt]{article}
\usepackage{amsmath,amsthm,amssymb}
\usepackage[matrix,arrow,curve]{xy}
\usepackage{verbatim}
\usepackage{mathtools}

\DeclarePairedDelimiter\ceil{\lceil}{\rceil}
\newtheorem{Lemma}{Lemma}[section]
\newtheorem{Def}[Lemma]{Definition}
\newtheorem{prop}[Lemma]{Proposition}
\newtheorem{theorem}[Lemma]{Theorem}
\newtheorem*{theorem*}{Theorem}
\newtheorem{coro}[Lemma]{Corollary}
\newtheorem{rema}[Lemma]{Remark}
\newcommand{\Hom}{{\rm Hom}}
\newcommand{\Id}{{\rm Id}}
\newcommand{\charr}{{\rm char\,}}

\DeclareMathOperator{\Ker}{Ker}

\DeclareMathOperator{\Aut}{Aut}
\renewcommand{\Im}{{\rm Im\,}}
\newcommand{\umod}{{\rm \underline{mod}}\mbox{-}}

\newcommand{\Mod}{{\rm Mod}\mbox{-}}
\newcommand{\Bimod}{{\rm Bimod}\mbox{-}}

\newcommand{\End}{{\rm End}}
\newcommand{\Nat}{{\rm Nat}}

\newcommand{\td}{\widetilde}
\newcommand{\ZZ}{{\rm Z}}
\newcommand{\Ab}{{\rm Ab}}

\newcommand{\HH}{{\rm HH}}
\newcommand{\HK}{{\rm HK}}
\newcommand{\HR}{{\rm HR}}

\newcommand{\bZ}{\mathbb{Z}}
\newcommand{\AAA}{{\rm A}}
\newcommand{\BB}{{\rm B}}

\newcommand{\DD}{{\rm D}}

\newcommand{\XX}{{\rm X}}
\newcommand{\YY}{{\rm Y}}

\newcommand{\NN}{{\rm N}}
\newcommand{\MM}{{\rm M}}

\newcommand{\SSS}{{\rm \Sigma}}

\newcommand{\UU}{{\rm U}}
\newcommand{\VV}{{\rm V}}
\newcommand{\WW}{{\rm W}}

\newcommand{\Kb}{{\rm K^b_p}}

\newcommand{\Ann}{{\rm Ann}}

\renewcommand{\le}{\leqslant}
\renewcommand{\ge}{\geqslant}

\newcommand{\ee}{\varepsilon}

\newcommand{\ph}{\varphi}
\newcommand{\ot}{\otimes}
\newcommand{\kk}{\mathbf{k}}
\newcommand{\ii}{\mathbf{i}}

\newenvironment{Proof}[1][Proof.]{\begin{trivlist}
\item[\hskip \labelsep {\bfseries #1}]}{\flushright
$\Box$\end{trivlist}}

\numberwithin{equation}{section}

\begin{document}
\title{K\"ulshammer ideals of graded categories and Hochschild cohomology.}
\author{Yury Volkov and Alexandra Zvonareva\thanks{Yury Volkov is supported by the RFBR Grant 17-01-00258 and by the President's Program "Support of Young Russian Scientists" (Grant MK-1378.2017.1). Alexandra Zvonareva is supported by the RBFR Grant 16-31-60089.}}
\date{}
\maketitle

\begin{abstract}
We generalize the notion of K\"ulshammer ideals to the setting of a graded category. This allows us to define and study some properties of K\"ulshammer type ideals in the graded center of a triangulated category and in the Hochschild cohomology of an algebra, providing new derived invariants. Further properties of K\"ulshammer ideals are studied in the case where the category is $d$-Calabi-Yau.  
\end{abstract}

\section{Introduction}
Let $\Lambda$ be a symmetric algebra over a field of positive characteristic $p$ with a symmetrizing form $(-,-)$. 
The sequence of K\"ulshammer ideals $$K^{cl}_r=\{a\in \Lambda | (a,b)=0 \mbox{ for } b\in\Lambda \mbox{ such that } b^{p^r}\in [\Lambda,\Lambda]\}$$ in the center of $\Lambda$ is a fine invariant of the derived category of an algebra \cite{Zim1, KLZ}. These ideals were applied to distinguish various algebras up to derived equivalence \cite{Hol1, Hol2}. With the use of trivial extensions the definition of K\"ulshammer ideals was extended to arbitrary algebras \cite{BZ}. Also various attempts to generalize K\"ulshammer ideals to higher Hochschild (co)homology were taken (see \cite{Zim2, Zim3}).

In this paper we propose to consider the same type of ideals in the center of a graded category $\AAA$. These ideals are defined using the module structure of the graded abelianization of the category $\AAA$ over its graded center. For a graded category $\AAA$ the center of $\AAA$ is the graded $\kk$-algebra $\AAA^{\AAA}$, whose $i$-th component $\AAA^{\AAA}_i$ is formed by elements $(m_{i,x})_{x\in\AAA}\in\Pi_{x\in\AAA}\AAA(x,x)_i$ such that $fm_{i,x}=(-1)^{ij}m_{i,y}f$ for any $x,y\in\AAA$ and any $f\in\AAA(x,y)_j$. The abelianization of $\AAA$ is the graded $\kk$-module $\AAA_{\AAA}=\oplus_{x\in\AAA}\AAA_x^x/[\AAA,\AAA]$, where $[\AAA,\AAA]$ denotes the subspace of $\oplus_{x\in\AAA}\AAA(x,x)$ formed by the elements  $fu-(-1)^{ij}uf$, for all $x,y\in\AAA$, $u\in\AAA(y,x)_i$ and $f\in\AAA(y,x)_j$.

The ideals $K_{r,s}\AAA$ are defined as the annihilators of the appropriate homogeneous component of the kernel of the map $(-)^{p^r}$. For precise definitions see Sections \ref{centab} and \ref{Kidealinc}. It turns out that the ideals defined in this way are invariant under graded equivalences.   

This construction is then applied in the particular situation of a category $\AAA$ with an automorphism $\SSS$. In this case the orbit category $\AAA/\SSS$ is graded and we can consider ideals in the graded center of $\AAA/\SSS$. If $\AAA$ is $d$-Calabi-Yau, we establish a duality between the Hochshild-Mitchel homology and cohomology of $\AAA$. This generalizes the well known duality between Hochschild homology and cohomology for symmetric algebras. Using the pairing provided by this duality we recover the usual definition of the ideals $K_{r,d}(\AAA/\SSS)$ in this general context. Thus,
$$K_{r,d}(\AAA/\SSS)=\{a\in (\AAA/\SSS)^{\AAA/\SSS} | (a,b)=0 \mbox{ for } b \in \oplus_{x\in \AAA}(\AAA/\SSS)_x^x \mbox{ such that } b^{p^r}\in [\AAA/\SSS,\AAA/\SSS]\}.$$

If $\AAA$ is the homotopy category of complexes of finitely generated projective modules $\Kb \Lambda$ over some algebra $\Lambda$ and $[1]$ is the shift functor, then $(\Kb \Lambda/[1])^{\Kb \Lambda/[1]}$ is the graded center $\ZZ^*(\Kb \Lambda)$ of $\Kb \Lambda$, i.e. the graded ring of natural transformations from the identity functor to $[n]$ which commute modulo the sign $(-1)^n$ with $[1]$. Graded centers of triangulated categories are studied and used by various authors and attract some interest (see \cite{BF, Li, BIK, KY}). Thus, the  ideals defined for $\Kb \Lambda/[1]$ belong to the graded center of $\Kb \Lambda$. 

The characteristic map from the Hochshild cohomology of $\Lambda$ to the graded center of $\Kb \Lambda$ allows us to define K\"ulshammer ideals $\HK_{r,s}(\Lambda)$ in the higher Hochshild cohomology as the inverse image of the ideals $K_{r,s}(\Kb \Lambda/[1])$. If $\Lambda$ is a  symmetric algebra,  then $\HK^0_{r,0}(\Lambda)$ coincide with the classical K\"ulshammer ideals $K^{cl}_r$. So the ideals $\HK_{r,s}(\Lambda)$ can be considered as a generalization of K\"ulshammer ideals to higher Hochshild cohomology. Both ideals $K_{r,s}(\Kb \Lambda/[1])$ and $\HK_{r,s}(\Lambda)$ are invariant under derived equivalence.

For arbitrary algebras we provide an alternative description of $\HK^0_{r,0}(\Lambda)$.
$$\HK^0_{r,0}(\Lambda)=\{a\in Z(\Lambda) | ab \in [\Lambda,\Lambda] \mbox{ for } b \mbox{ such that } b^{p^r}\in [\Lambda,\Lambda]\}.$$
This provides an alternative generalization of K\"ulshammer ideals for non-symmetric algebras. 

We finish the paper computing all the defined ideals for the algebra $k[x]/x^2$.

\section{Graded center and abelianization}\label{centab}

All categories and functors are assumed to be $\kk$-linear for some fixed field $\kk$. Moreover, all categories are assumed to be small. We write simply $\ot$ and $\Hom$ instead of $\ot_{\kk}$ and $\Hom_{\kk}$.
In this section we recall some definitions and introduce some notation.
From here on we will write $\AAA_x^y$ instead of $\Hom_{\AAA}(x,y)$ for $\Hom$-sets in a category $\AAA$.

\begin{Def}{\rm A category $\AAA$ is called a {\it graded category} if, for any $x$, $y$ in $\AAA$, there is a fixed decomposition of graded spaces $\AAA_x^y=\oplus_{i\in\bZ}(\AAA_x^y)_i$ such that $fg\in(\AAA_x^z)_{i+j}$ for any $f\in(\AAA_y^z)_j$ and $g\in(\AAA_x^y)_i$. We write $|f|$ for the degree of $f\in\AAA_x^y$, i.e. $|f|=i$ if and only if $f\in(\AAA_x^y)_i$.
}
\end{Def}

\begin{Def} {\rm The {\it tensor product} $\AAA\ot\BB$ of categories $\AAA$ and $\BB$ is a
$\kk$-linear category defined in the following way. Its objects are
pairs $(x,y)$ where $x\in\AAA$, $y\in\BB$. Its morphism spaces are
$$
(\AAA\ot\BB)_{(x_1,y_1)}^{(x_2,y_2)}=\AAA_{x_1}^{x_2}\ot\BB_{y_1}^{y_2}.
$$
The composition in $\AAA\ot\BB$ is given by the formula
$$
(f_2\ot g_2)(f_1\ot g_1)=f_2f_1\ot g_2g_1,
$$
where $f_1\in\AAA_{x_1}^{x_2}$, $f_2\in\AAA_{x_2}^{x_3}$,
$g_1\in\BB_{y_1}^{y_2}$, $g_2\in\BB_{y_2}^{y_3}$ and $x_1,x_2,x_3\in\AAA$,
$y_1,y_2,y_3\in\BB$.}
\end{Def}

\begin{Def}{\rm $\kk$-linear contravariant functors from $\AAA$ to $\Mod\kk$ are called {\it
$\AAA$-modules}. We denote by $\Mod\AAA$ the category of
$\AAA$-modules.
An {\it $\AAA$-bimodule} is by definition an $\AAA^{\rm op}\ot\AAA$-module. We denote by $\Bimod\AAA$ the category of $\AAA$-bimodules.
}
\end{Def}

If it does not cause any confusion, we will write $fug$ instead of $\MM(f\otimes g)u$ for an $\AAA$-bimodule $\MM$ and $g\in\AAA_z^w$, $f\in\AAA_x^y$, $u\in\MM(x,w)$, where $w,x,y,z\in\AAA$.

If $\AAA$ is graded, then the $\AAA$-bimodule $\MM$ is called {\it graded} if there are some fixed decompositions $\MM(x,y)=\oplus_{i\in\bZ}\MM(x,y)_i$ such that $fug\in\MM(y,z)_{i+j+l}$ for all $g\in(\AAA_z^w)_j$, $f\in(\AAA_x^y)_i$ and $u\in\MM(x,w)_l$.
A homogeneous morphism of degree $m$ from a graded bimodule $\MM$ to a graded bimodule $\NN$ is a collection of maps $\phi_{x,y,k}:\MM(x,y)_k\rightarrow \NN(x,y)_{k+m}$ such that $\NN(f\otimes g)\phi_{x,w,k}=(-1)^{mi}\phi_{y,z,k+i+j}\MM(f\otimes g)$ for $g\in(\AAA_z^w)_j$, $f\in(\AAA_x^y)_i$. A morphism between graded bimodules is by definition a finite sum of homogeneous morphisms. 
Any nongraded category $\AAA$ can be considered as a graded category with all morphisms of degree 0. In this case a graded bimodule is simply a module $\MM$ with $\AAA$-bimodule decomposition $\MM=\oplus_{i\in\bZ}\MM_i$.
We define the graded bimodule $\MM[n]$ as a graded $\AAA$-bimodule such that $\MM[n](x,y)_i=\MM(x,y)_{i+n}$ and $\MM[n](f\otimes g)u=(-1)^{in}\MM(f\otimes g)u$ for all $g\in(\AAA_z^w)_j$, $f\in(\AAA_x^y)_i$ and $u\in\MM(x,w)_l$. If $\phi_{x,y,j}$ is a homogeneous morphism form $\MM$ to $\NN$, then $\phi_{x,y,j}[n]:=\phi_{x,y,j+n}$.

One can consider  $\AAA$ as an $\AAA$-bimodule defined by the equality $\AAA(x,y)=\AAA_y^x$ on objects and in the obvious way on morphisms.


\begin{Def}\label{lincat}
{\rm An {\it $\AAA$-linear category} is an $\AAA$-bimodule $\MM$ together with a structure of a $\kk$-linear category (which will be also denoted by $\MM$) compatible with the bimodule structure. Namely, the class of objects of $\MM$ is the same as the class of objects of $\AAA$, and the morphism spaces of $\MM$ are $\MM_y^x=\MM(x,y)$. Thus, there are bilinear maps $$\MM(y,x)\times \MM(x,z)= \MM_x^y \times \MM_z^x \xrightarrow{-\circ -} \MM_z^y= \MM(y,z)$$ that satisfy all the conditions of a categorical composition. The compatibility conditions are
$$
f(v \circ u)g= (fv)\circ (ug)\mbox{ and }u \circ (fv)= (uf) \circ v,
$$
$$
\xymatrix{ \MM(y,x)\times \MM(x,z) \ar[d]^{\MM(f\otimes id)\times \MM(id\otimes g)} \ar[r]^-{-\circ-} &\MM(y,z)\ar[d]^{\MM(f\otimes g)}&&\MM(z,w)\times \MM(y,x) \ar[d]^{id\times\MM(f\otimes id)} \ar[rr]^-{\MM(id\otimes f)\times id} &&\MM(z,y)\times \MM(y,x)\ar[d]^-{-\circ-}\\
 \MM(w,x)\times \MM(x,\tau) \ar[r]^-{-\circ-} &\MM(w,\tau)&&\MM(z,w)\times \MM(w,x)  \ar[rr]^-{-\circ-} &&\MM(z,y)\times \MM(y,x)}
 $$
where $f \in \AAA_y^w $ and $g \in \AAA_{\tau}^z$ are morphisms in $\AAA$ and $u$ and $v$ are morphisms in $\MM$.
}
\end{Def}

We say that an $\AAA$-linear category $\MM$ is graded if $\MM$ is a graded category and a graded bimodule with respect to the same family of decompositions $\MM^x_y=\oplus_{i\in\bZ}(\MM^x_y)_i$.

\begin{Def}{\rm Let $\AAA$ be a graded category and let $\MM$ be a graded $\AAA$-bimodule. The {\it $\AAA$-center} of $\MM$ is the graded $\kk$-module $\MM^{\AAA}$, whose $i$-th component $\MM^{\AAA}_i$ is formed by such elements $(m_{i,x})_{x\in\AAA}\in\Pi_{x\in\AAA}\MM(x,x)_i$ that $fm_{i,x}=(-1)^{ij}m_{i,y}f$ for any $x,y\in\AAA$ and any $f\in(\AAA_x^y)_j$.
If $\AAA$ and $\MM$ are not graded, then we can consider them as a graded category and a bimodule concentrated in degree 0. So we can talk about the center of a nongraded bimodule over a nongraded category.
}
\end{Def}

Note that Definition \ref{lincat} guarantees that if $\MM$ is an $\AAA$-linear category, then $\MM^{\AAA}$ has a structure of a unital associative $\kk$-algebra, with multiplication induced by the composition in $\MM$. In particular, the compatibility of the the bimodule and the categorical structure ensures that $\MM^{\AAA}$ is closed under multiplication.

\begin{Def}{ \rm Let $\AAA$ be a (nongraded) category and $\SSS:\AAA\rightarrow \AAA$ be an automorphism of $\AAA$. The {\it orbit category} $\AAA/\SSS$ is a graded category defined as follows.
\begin{itemize}
\item The class of objects of $\AAA/\SSS$ is equal to that of $\AAA$;
\item The sets of morphisms are  $\big((\AAA/\SSS)_x^y\big)_n=\AAA_x^{\SSS^ny}$ for $x,y\in\AAA/\SSS$ and $n\in\bZ$;
\item The composition $\circ$ in $\AAA/\SSS$ is given by the formula
$g\circ f=\SSS^n(g) f$ for $f\in\big((\AAA/\SSS)_x^y\big)_n$ and $g\in\big((\AAA/\SSS)_y^z\big)_m$.
\end{itemize}}
\end{Def}

Note that $\AAA/\SSS$ becomes a graded $\AAA$-linear category if we define $(\AAA/\SSS)(f\otimes g)u=\SSS^n(f)ug$ for $f\in\AAA_{y_1}^{y_2}$, $g\in\AAA_{x_1}^{x_2}$ and $u\in\big((\AAA/\SSS)_{x_2}^{y_1}\big)_n$. Given an automorphism $\alpha$ of the category $\AAA$ and $\MM\in\Bimod\AAA$, we define the bimodule ${}^{\alpha}\MM$ as the composition of functors $\MM\circ (\alpha\ot\Id_{\AAA})\in\Bimod\AAA$. It is easy to see that this defines an action of $(\Aut\AAA)^{\rm op}$ on $\Bimod\AAA$. 
Note that $\AAA/\SSS\cong\oplus_{n\in\mathbb{Z}}\big({}^{\SSS^n}\AAA\big)[-n]$ as a graded  $\AAA$-bimodule. 

If $\alpha$ is an automorphism of $\AAA$, then we say that $\alpha$ acts on the graded bimodule $\MM$ if there is a homogeneous isomorphism $ \MM \xrightarrow{\alpha_{\MM}} \MM \circ(\alpha \ot \alpha)$ of degree $0$. Such action induces an automorphism $\alpha_{\MM^{\AAA}}$ of the graded center $\MM^{\AAA}$. If $\MM$ is an $\AAA$-linear category and $\alpha$ acts on it by a category automorphism, then $\alpha_{\MM^{\AAA}}$ is an automorphism of the graded algebra $\MM^{\AAA}$. Note that if $\alpha$ acts on $\MM$ by $\alpha_{\MM}$, then we can define $\alpha_{\MM[n]}$ as $\alpha_{\MM[n]}=(-1)^n\alpha_{\MM}[n]$.
The automorphism $\SSS$ of $\AAA$ acts on $^{\SSS^n}\AAA$ by the rule $\SSS_{(^{\SSS^n}\AAA)}(u)=\SSS(u)$ for $u\in \AAA_x^{\SSS^ny}$. We fix this action of $\SSS$ on $^{\SSS^n}\AAA$ and  write simply $\SSS$ instead of $\SSS_{(^{\SSS^n}\AAA)}$.
Moreover, for any integer $m$, we denote by $\SSS$ the natural transformation $\SSS_{({}^{\SSS^n}\AAA)[m]}$ and the automorphism $\SSS_{\big(({}^{\SSS^n}\AAA)[m]\big)^{\AAA}}$.
Thus, $\SSS$ acts on the bimodule $\AAA/\SSS$ and this action determines an automorphism of the graded $\AAA$-linear category $\AAA/\SSS$. For any space $\VV$ with an action of some automorphism $\alpha$ we can consider the subspace of invariants $\VV^{\alpha}$, i.e $\{m\in \VV | \alpha(m)=m\}$.

\begin{Def}{\rm
Let $\AAA$ be a (nongraded) category and let $\SSS:\AAA\rightarrow \AAA$ be some fixed automorphism of $\AAA$. We define the graded rings $\Nat^*(\AAA)=\Nat^*(\AAA,\SSS)$ and $\ZZ^*(\AAA)=\ZZ^*(\AAA,\SSS)$ as follows. Let $\Nat^n(\AAA)$ ($n\in\mathbb{Z}$)
be the abelian group formed by all natural transformations $\eta: \Id_{\AAA}\rightarrow\SSS^n$ and 
 $\ZZ^n(\AAA)$ be its subgroup formed by $\eta$ that satisfy the equality
$\eta\SSS = (-1)^n\SSS\eta$. Given natural transformations $\eta:\Id_{\AAA} \rightarrow\SSS^n$ and $\theta:\Id_{\AAA} \rightarrow\SSS^m$ ($n,m\in\mathbb{Z}$ ), we define the product of $\eta$ and $\theta$ by the formula $\eta\theta=\SSS^m(\eta)\circ\theta:\Id_{\AAA}\rightarrow\SSS^{n+m}$. We call $\ZZ^*(\AAA)$ the {\it graded center} of $\AAA$. Note that $\eta\theta = (-1)^{mn}\theta\eta$ if $\eta\in \ZZ^*(\AAA)$.
}
\end{Def}
\begin{rema} The definition of a graded center has sense if $\SSS$ is an autoequivalence, but further we need it to be an automorphism. On the other hand, we can replace any autoequivalence of a category by an automorphism (the category is changed during this process) by the results of \cite{Asashiba}.
\end{rema}

\begin{Lemma}\label{gc}
There is an isomorphism of graded algebras
$\Nat^*(\AAA)\cong (\AAA/\SSS)^{\AAA}$ that induces isomorphisms $\ZZ^*(\AAA)\cong\Big( (\AAA/\SSS)^{\AAA}\Big)^{\SSS}\cong(\AAA/\SSS)^{\AAA/\SSS}$.
\end{Lemma}
\begin{Proof}
By definition, $\AAA$-center of $\AAA/\SSS$ has as the $n$-th component $\big((\AAA/\SSS)^{\AAA}\big)_n$  elements $(m_{n,x})_{x\in\AAA}\in\Pi_{x\in\AAA}(\AAA/\SSS)(x,x)_n$ such that $\SSS^n(f)m_{n,x}=m_{n,y}f$
for any $x,y\in\AAA$ and any $f\in\AAA_x^y$.
 Thus, a family $(m_{n,x})_{x\in\AAA}$ belongs to $\big((\AAA/\SSS)^{\AAA}\big)_n$ if and only if it gives a natural transformation $\Id_{\AAA} \rightarrow \SSS^n$. The multiplication on $\Nat^*(\AAA)$ is $\eta\theta=\SSS^m(\eta)\circ\theta$, the multiplication on $(\AAA/\SSS)^{\AAA}$ is induced by the composition in $\AAA/\SSS$, i.e. $gf=\SSS^n(g)\circ f$, hence these graded algebras are isomorphic.

As for the isomorphism $\ZZ^*(\AAA)\cong(\AAA/\SSS)^{\AAA/\SSS}$, it is clear that an element of $\big((\AAA/\SSS)^{\AAA/\SSS}\big)_n$ is a natural transformation, as in the previous paragraph, i.e. $(\AAA/\SSS)^{\AAA/\SSS}$ is a subalgebra of $(\AAA/\SSS)^{\AAA}$. Taking $y=\SSS x$ and $f=\Id_x\in\AAA_x^{\SSS^{-1}y}=\big((\AAA/\SSS)_x^y\big)_{-1}$ we get $$m_{n,x}=(-1)^n \SSS^{-1}( m_{n,y})=(-1)^n \SSS^{-1}(m_{n,\SSS x})$$ and hence $\SSS(m_{n,x})=(-1)^n m_{n,\SSS x}$; thus, we get an element of $\ZZ^*(\AAA)$. If we take an element $\eta$ of $\ZZ^n(\AAA)$, then for $f: x \rightarrow \SSS^k y$ we get $\SSS^n(f)\eta_x =\eta_{\SSS^k y}f=(-1)^{nk}\SSS^k(\eta_y)f$; thus, we get an element of $\big((\AAA/\SSS)^{\AAA/\SSS}\big)_n$.

The isomorphism $\ZZ^*(\AAA)\cong((\AAA/\SSS)^{\AAA})^{\SSS}$ follows from the definition of the action of $\SSS$ on the bimodule $\AAA/\SSS$: the family $\{m_{n,x}\}_{x\in\AAA}$ belongs to $((\AAA/\SSS)^{\AAA})^{\SSS}$ if and only if $(-1)^n\SSS(m_{m,x})= m_{n,\SSS x}$ if and only if the corresponding natural transformation belongs to $\ZZ^*(\AAA)$.
\end{Proof}

\begin{Def}{\rm
Given a graded $\AAA$-bimodule $\MM$, $[\AAA,\MM]$ denotes the subspace of $\oplus_{x\in\AAA}\MM_x^x$ formed by the elements  $fu-(-1)^{ij}uf$, for all $x,y\in\AAA$, $u\in\MM(y,x)_i$ and $f\in(\AAA_y^x)_j$. The {\it $\AAA$-abelianization} of $\MM$ is the graded $\kk$-module $\MM_{\AAA}=\oplus_{x\in\AAA}\MM_x^x/[\AAA,\MM]$. As in the case of the center, we can talk about the abelianization of a nongraded bimodule over a nongraded category.
}
\end{Def}

For any space $\VV$ with an action of some automorphism $\alpha$ we can consider the space of co-invariants $\VV_{\alpha}$, i.e. the quotient space of $\VV$ modulo the subspace generated by the classes of the elements of the form $v-\alpha(v)$ for $v\in\VV$. If $\MM$ is an $\AAA$-linear category, then the composition in $\MM$ induces a structure of a graded $\MM^{\AAA}$-bimodule on $\MM_{\AAA}$. Note that an action of an automorphism on $\MM$ induces an action on $\MM_{\AAA}$. So in this case one can define the graded space
$(\MM_{\AAA})_{\SSS}$ as the quotient space of $\MM_{\AAA}$ modulo the subspace generated by the classes of the elements of the form $m-\SSS m$ for $m\in\oplus_{x\in\AAA}\MM_x^x$. The $\MM^{\AAA}$-bimodule structure on $\MM_{\AAA}$ induces the 
$(\MM^{\AAA})^{\SSS}$-bimodule structure on $(\MM_{\AAA})_{\SSS}$. Note also that if the action of the automorphism $\alpha$ of $\AAA$ on the bimodule $\MM$ is given by $\alpha_{\MM}$, then $\alpha$ acts on $\MM$ by $a\alpha_{\MM}$  for any $a\in\kk^*$ as well.

\begin{Def}{\rm
Let $\AAA$ be a (nongraded) category and let $\SSS:\AAA\rightarrow \AAA$ be some fixed automorphism of $\AAA$. We define the graded spaces $\Nat_*(\AAA)=\Nat_*(\AAA,\SSS)$ and $\Ab_*(\AAA)=\Ab_*(\AAA,\SSS)$ as follows.
Let $\Nat_n(\AAA)$ be the quotient space of the space $\oplus_{x\in\AAA}\AAA_x^{\SSS^nx}$ modulo the subspace generated by the elements of the form $fg-\SSS^n(g)f$ for all $f\in\AAA_x^{\SSS^{n}y}$ and $g\in\AAA_y^{x}$.
We define $\Ab_n(\AAA)$ as the quotient space of $\Nat_n(\AAA)$ modulo the subspace generated by the classes of the elements of the form $f+(-1)^n\SSS(f)$ for all $f\in\AAA_x^{\SSS^nx}$. We call $\Ab_*(\AAA)$ the {\it graded abelianization} of $\AAA$.
}
\end{Def}

It is easy to see that $\Nat_*(\AAA)$ is a $\Nat^*(\AAA)$-bimodule. Moreover, the corresponding $\ZZ^*(\AAA)$-bimodule structure on $\Nat_*(\AAA)$ induces a $\ZZ^*(\AAA)$-bimodule structure on $\Ab_*(\AAA)$.
Note also that the isomorphisms from Lemma \ref{gc} induce a $\Nat^*(\AAA)$-bimodule structure on $(\AAA/\SSS)_{\AAA}$ and a $\ZZ^*(\AAA)$-bimodule structure on $(\AAA/\SSS)_{\AAA/\SSS}$ and $\Big( (\AAA/\SSS)_{\AAA}\Big)_{-\SSS}$ (in fact, as it was noted above, on $\Big( (\AAA/\SSS)_{\AAA}\Big)_{a\SSS}$ for any $a\in\kk^*$).

\begin{Lemma}\label{gab}
There is an isomorphism of graded $\Nat^*(\AAA)$-bimodules
$\Nat_*(\AAA)\cong (\AAA/\SSS)_{\AAA}$ that induces isomorphisms $\Ab_*(\AAA)\cong\Big( (\AAA/\SSS)_{\AAA}\Big)_{-\SSS}\cong(\AAA/\SSS)_{\AAA/\SSS}$.
\end{Lemma}
\begin{Proof} By definition, $\Nat_*(\AAA)$ and $(\AAA/\SSS)_{\AAA}$ are both quotients of the $\Nat^*(\AAA)$-bimodule $\oplus_{n\in\mathbb{Z}}\oplus_{x\in\AAA}\AAA_x^{\SSS^nx}$ and it suffices to check that we take the quotient space modulo the same submodule.

It is easy to see that  $$fg-\SSS^n(g)f=(\AAA/\SSS)(\Id_y\ot g)(f)-(\AAA/\SSS)(g\ot \Id_x)(f)$$ for all $f\in\AAA_x^{\SSS^{n}y}$ and $g\in\AAA_y^{x}$, i.e. $\Nat_n(\AAA)=\big((\AAA/\SSS)_{\AAA}\big)_n=\oplus_{x\in\AAA}\AAA_x^{\SSS^nx}/\UU_n$, where $\UU_n$ is generated by $fg-\SSS^n(g)f$ for  $f\in\AAA_x^{\SSS^{n}y}$ and $g\in\AAA_y^{x}$.

Now we have $\Ab_n(\AAA)=\oplus_{x\in\AAA}\AAA_x^{\SSS^nx}/(\UU_n+\VV_n)$, where $\VV_n$ is generated by $f+(-1)^n\SSS(f)$ for all $f\in\AAA_x^{\SSS^nx}$. It is easy to see that $\Big(\big( (\AAA/\SSS)_{\AAA}\big)_{-\SSS}\Big)_n=\oplus_{x\in\AAA}\AAA_x^{\SSS^nx}/(\UU_n+\VV_n)$ too. Let now consider $\big((\AAA/\SSS)_{\AAA/\SSS}\big)_n=\oplus_{x\in\AAA}\AAA_x^{\SSS^nx}/\WW_n$, where $\WW_n$ is generated by $\SSS^j(f)g-(-1)^{ij}\SSS^i(g)f$ for  $f\in\AAA_x^{\SSS^{i}y}$ and $g\in\AAA_y^{\SSS^{j}x}$, $i+j=n$. Taking $j=0$, we get that $\UU_n\subseteq\WW_n$. Taking $i=n+1$, $j=-1$, $y=\SSS^{-1}x$, and $g=\Id_x$, we get that  $\VV_n\subseteq\WW_n$. Thus, $\UU_n+\VV_n\subseteq\WW_n$.

Since $\SSS^j(f)g-\SSS^n(g)\SSS^j(f)\in\UU_n$ and
$$
\SSS^i(g)f-(-1)^{ij}\SSS^n(g)\SSS^j(f)=\SSS^i(g)f-(-1)^{j(n+1)}\SSS^j(\SSS^i(g)f)\in\VV_n,
$$
we have $\WW_n=\UU_n+\VV_n$.
\end{Proof}

\section{K\"ulshammer ideals in the center of a graded category}\label{Kidealinc}

In this section we define K\"ulshammer ideals in the center of a graded category. From now on we assume that $\kk$ is a field of characteristic $p>0$.
For a graded category $\AAA$ we define the map $$\xi_p:\oplus_{x\in\AAA}\AAA_x^x\rightarrow\oplus_{x\in\AAA}\AAA_x^x$$ by the equality $\xi_p\left(\sum_{n\in\mathbb{Z}}f_n\right)=\sum_{n\in\mathbb{Z}}f_n^p$, where $f_n\in\left(\oplus_{x\in\AAA}\AAA_x^x\right)_n$ for $n\in\mathbb{Z}$.
Note that $\xi_p$ maps $(\AAA_x^x)_n$ to $(\AAA_x^x)_{np}$.

\begin{Lemma}\label{ppower}
Let $\AAA$ be a graded category. The map $\xi_p:\oplus_{x\in\AAA}\AAA_x^x\rightarrow \oplus_{x\in\AAA}\AAA_x^x$ induces a well-define map $\xi_p:\AAA_{\AAA}\rightarrow \AAA_{\AAA}$.
\end{Lemma}
\begin{Proof} It is enough to show that $\xi_p\left([\AAA,\AAA]_n\right)\subseteq [\AAA,\AAA]_{np}$. Let us first prove that for $f,g\in \left(\oplus_{x\in\AAA}\AAA_x^x\right)_n$ one has $(f+g)^p-f^p-g^p\in [\AAA,\AAA]_{np}$.
Fix $s_i\in\left(\oplus_{x\in\AAA}\AAA_x^x\right)_n$ for $1\le i\le p$. By definition of $[\AAA,\AAA]$, we have $s_1\dots s_p-(-1)^{n^2(p-1)}s_2\dots s_ps_1\in [\AAA,\AAA]$. Since $(-1)^{n^2(p-1)}=1$ if $2\nmid p$ or $\charr\kk=p=2$, we have $s_1\dots s_p-s_2\dots s_ps_1\in [\AAA,\AAA]$ for any prime $p$.
Now the condition $(f+g)^p-f^p-g^p\in [\AAA,\AAA]_{np}$ can be proved in the same way as in \cite[Lemma 1.1]{Zim}.

Any element of $[\AAA,\AAA]_n$ has the form $u=\sum_{t=1}^k\big(f_tg_t-(-1)^{i_t(n-i_t)}g_tf_t\big)$ for some integer $k\ge 0$ and some $f_t\in \big(\AAA_{x_t}^{y_t}\big)_{i_t}$ and $g_t\in \big(\AAA_{y_t}^{x_t}\big)_{n-i_t}$ for $1\le t\le k$. Then we have
\begin{multline*}
\xi_p(u)+[\AAA,\AAA]_{np}=\sum_{t=1}^k\big((f_tg_t)^p-(-1)^{i_t(n-i_t)}(g_tf_t)^p\big)+[\AAA,\AAA]_{np}\\
=\sum_{t=1}^k\Big(f_t\big((g_tf_t)^{p-1}g_t\big)-(-1)^{i_t(np-i_t)}\big((g_tf_t)^{p-1}g_t\big)f_t\Big)+[\AAA,\AAA]_{np}\in [\AAA,\AAA]_{np},
\end{multline*}
i.e. $\xi_p\left([\AAA,\AAA]_n\right)\subseteq [\AAA,\AAA]_{np}$.
\end{Proof}

\begin{Def}{\rm
Let $\AAA$ be a graded category. Then we define $T_r\AAA\subseteq \AAA_{\AAA}$ as the kernel of the map $\xi_p^r$.
Let us define the graded subspace $K_{r,s}\AAA\subseteq \AAA^{\AAA}$ by the equality $$(K_{r,s}\AAA)_n=\Ann_{\AAA^{\AAA}}\big((T_r\AAA)_{s-n}\big)_n.$$
Here, for an algebra $\Lambda$, a $\Lambda$-bimodule $\MM$ and $\VV\subseteq \MM$, we denote by $\Ann_{\Lambda}(\VV)$ the set $\{a\in \Lambda\mid aV=Va=0\}$. 
We will call $K_{r,s}\AAA$ the $(r,s)$-th {\it K\"ulshammer ideal} of $\AAA$ and $K_r\AAA=\cap_{s\in \mathbb{Z}}K_{r,s}\AAA=\Ann_{\AAA^{\AAA}}(T_r\AAA)$ the $r$-th {\it K\"ulshammer ideal} of $\AAA$.
Also we will call $R_s\AAA=\cap_{r\ge 0}K_{r,s}(\AAA)$ the $s$-th {\it Reynolds ideal} of $\AAA$ and will call $R\AAA=\cap_{r\ge 0}K_r(\AAA)$ the {\it Reynolds ideal} of $\AAA$.
}
\end{Def}

We say that $F:\AAA\rightarrow\AAA'$ is a {\it degree preserving equivalence} if $F$ induces an isomorphism of graded spaces $F:\AAA_x^y\rightarrow (\AAA')_{Fx}^{Fy}$ for any $x,y\in\AAA$, and for any $x'\in\AAA'$ there exists an isomorphism $\xi_{x'}\in\AAA_{x'}^{Fx}$ of degree zero for some $x\in\AAA$.

\begin{theorem}\label{ideals}
If $\AAA$ is a graded category and $s\in\mathbb{Z}$, then
$$
\AAA^{\AAA}= K_{0,s}\AAA\supseteq K_{1,s}\AAA\supseteq\dots\supseteq K_{r,s}\AAA\supseteq\dots\supseteq R_s\AAA
$$
is a decreasing sequence of graded ideals. In particular,
$$
\AAA^{\AAA}= K_0\AAA\supseteq K_1\AAA\supseteq\dots\supseteq K_r\AAA\supseteq\dots\supseteq R\AAA
$$
is a decreasing sequence of graded ideals.
Moreover, if $F:\AAA\rightarrow\AAA'$ is a degree preserving equivalence, then $F$ induces an isomorphism of graded algebras $\ph_F:\AAA^{\AAA}\rightarrow (\AAA')^{\AAA'}$ such that $\ph_F(K_{r,s}\AAA)=K_{r,s}\AAA'$ for any $r\ge 0, s \in \mathbb{Z}$.
\end{theorem}
\begin{Proof} It follows directly from the definition that $R_s\AAA\subseteq K_{r+1,s}\AAA\subseteq K_{r,s}\AAA$ for all $r\ge 0$. Since $\AAA^{\AAA}$ is graded commutative, it is easy to see that
$\oplus_{n\in\mathbb{Z}}\big(\AAA_n^{\AAA}\cap \Ann_{\AAA^{\AAA}}(\VV_{s-n})\big)$ is a graded ideal of $\AAA^{\AAA}$ for any graded $\AAA^{\AAA}$-bimodule $\MM$ and any graded subbimodule $\VV\subseteq \MM$. For $f\in (\AAA_{\AAA})_m$ such that $f^{p^r}=0$, $\theta\in (\AAA^{\AAA})_{d-m-n}$ we have $(\theta f)^{p^r}=(-1)^{a}\theta^{p^r}f^{p^r}=0$, where $a=\frac{p^r(p^r-1)}{2}m(d-m-n)$, and hence $T_r\AAA$ is a graded subbimodule of $\AAA_{\AAA}$. Thus, the first part of the theorem is proved.

Let us now prove the second part. We can choose a quasi inverse equivalence $F'$ for $F$ and natural isomorphisms $\alpha:\Id_{\AAA}\rightarrow  F'F$ and $\beta:\Id_{\AAA'}\rightarrow FF'$ in such a way that $\alpha_x$ and $\beta_{x'}$ are of degree zero for all $x\in\AAA$ and $x'\in\AAA'$.

Suppose that $f=(f_{x})_{x\in\AAA}\in\Pi_{x\in\AAA}(\AAA_x^x)_n$ belongs to $\AAA^{\AAA}_n$. We define $\ph_F(f)_{x'}=\beta_{x'}^{-1}F(f_{F'x'})\beta_{x'}$. It is easy to see that $\ph_F(f)\in (\AAA')^{\AAA'}_n$. It is clear also that $\ph_F:\AAA^{\AAA}\rightarrow (\AAA')^{\AAA'}$ is a homomorphism of graded algebras.
Let us define $\ph_{F'}:(\AAA')^{\AAA'}\rightarrow \AAA^{\AAA}$ in the same way (using $\alpha$ instead of $\beta$). Take $f=(f_{x})_{x\in\AAA}\in\AAA^{\AAA}_n$. We have $f_{F'Fx}\alpha_x=\alpha_xf_x$ and $F'F(f_{x})\alpha_x=\alpha_xf_x$ since $f$ belongs to the center and $\alpha$ is a natural transformation. Hence, $f_{F'Fx}=F'F(f_x)$. Thus,
\begin{multline*}
\ph_{F'}\ph_F(f)_x=\alpha_x^{-1}F'\big(\beta_{Fx}^{-1}F(f_{F'Fx})\beta_{Fx}\big)\alpha_x\\
=\alpha_x^{-1}F'\big(\beta_{Fx}^{-1}FF'F(f_{x})\beta_{Fx}\big)\alpha_x=\alpha_x^{-1}F'F(f_{x})\alpha_x=f_x.
\end{multline*}
Thus, $\ph_{F'}\ph_F=\Id_{\AAA^{\AAA}}$. Analogously, $\ph_{F}\ph_{F'}=\Id_{(\AAA')^{\AAA'}}$. Consequently, $\ph_F$ is an isomorphism.

Suppose now that $f=(f_{x})_{x\in\AAA}\in\Pi_{x\in\AAA}(\AAA_x^x)_n$ belongs to $(K_{r,s}\AAA)_n$. Let us take some $u=(u_{x'})_{x'\in\AAA'}\in\big(\oplus_{x'\in\AAA'}(\AAA')_{x'}^{x'}\big)_{s-n}$ such that $u^{p^r}\in [\AAA',\AAA']$.
Let $(F')^{-1}(x)$ be the inverse image of $x$, i.e. the set $\{x'\in\AAA'\mid F'(x')=x\}$. We define $F'(u)_x=\sum\limits_{x'\in(F')^{-1}(x)}F'(u_{x'})$. By Lemma \ref{ppower}, we have $\big(F'(u)^{p^r}\big)_x=\sum\limits_{x'\in(F')^{-1}(x)}F'(u_{x'})^{p^r}=\sum\limits_{x'\in(F')^{-1}(x)}F'\big(u_{x'}^{p^r}\big)$.
Since $F'([\AAA',\AAA']) \subseteq [\AAA,\AAA]$, we see that $F'(u)\in T_r\AAA$. Then $fF'(u)\in [\AAA,\AAA]$. As before, we have $F(f)FF'(u)\in [\AAA',\AAA']$, where $\big(F(f)FF'(u)\big)_{x'}=\sum\limits_{FF'y'=x'}F(f_{F'y'})FF'(u_{y'})$. Since
$$\sum_{y'\in\AAA'}\big(F(f_{F'y'})FF'(u_{y'})-\beta_{y'}^{-1}F(f_{F'y'})FF'(u_{y'})\beta_{y'}\big)\in[\AAA',\AAA']$$
and $FF'(u_{y'})\beta_{y'}=\beta_{y'}u_{y'}$, we have
$$
\ph_F(f)u=F(f)FF'(u)-\sum_{y'\in\AAA'}\big(F(f_{F'y'})FF'(u_{y'})-\beta_{y'}^{-1}F(f_{F'y'})\beta_{y'}u_{y'}\big)\in[\AAA',\AAA'].
$$
Consequently, $\ph_F(f)\in (K_{r,s}\AAA')_n$, i.e. $\ph_F(K_{r,s}\AAA)\subseteq K_{r,s}\AAA'$. In the same way one can prove that $\ph_{F}^{-1}(K_{r,s}\AAA')=\ph_{F'}(K_{r,s}\AAA')\subseteq K_{r,s}\AAA$, i.e. $\ph_F(K_{r,s}\AAA)=K_{r,s}\AAA'$.
\end{Proof}

\begin{rema}
Note that the map $\ph_F$ constructed in the proof does not depend on $\beta$ and $F'$, i.e. it is really induced by $F$. Indeed, if $\beta,\td\beta:\Id_{\AAA'}\rightarrow FF'$ are two natural isomorphisms, then $\td\beta\beta^{-1}$ is a natural isomorphism from $\Id_{\AAA}$ to itself and hence 
$$(\td\beta\beta^{-1})_{x'}F(f_{F'x'})\beta_{x'}=F(f_{F'x'})\td\beta_{x'}\beta_{x'}^{-1}\beta_{x'}=F(f_{F'x'})\td\beta_{x'},$$ i.e. $\beta_{x'}^{-1}F(f_{F'x'})\beta_{x'}=\td\beta_{x'}^{-1}F(f_{F'x'})\td\beta_{x'}$. Now, if $F''$ is another quasi inverse of $F$, then there is a natural isomorphism $\gamma:F'\rightarrow F''$. Then $F(\gamma)\beta:\Id_{\AAA'}\rightarrow FF''$ is a natural isomorphism too and we have 
$$\big(F(\gamma_{x'})\beta_{x'}\big)^{-1}F(f_{F'x'})F(\gamma_{x'})\beta_{x'}=\beta_{x'}^{-1}F(\gamma_{x'}^{-1}f_{F'x'}\gamma_{x'})\beta_{x'}=\beta_{x'}^{-1}F(f_{F''x'})\beta_{x'}$$
since $f\in \AAA^{\AAA}$.
\end{rema}

\begin{coro}\label{ideals_aut}
If $\AAA$ is a category with the automorphism $\SSS$ and $s\in\mathbb{Z}$, then
$$
\ZZ^*(\AAA)= K_{0,s}(\AAA/\SSS)\supseteq K_{1,s}(\AAA/\SSS)\supseteq\dots\supseteq K_{r,s}(\AAA/\SSS)\supseteq\dots\supseteq R_s(\AAA/\SSS)
$$
is a decreasing sequence of graded ideals. In particular,
$$
\ZZ^*(\AAA)= K_0(\AAA/\SSS)\supseteq K_1(\AAA/\SSS)\supseteq\dots\supseteq K_r(\AAA/\SSS)\supseteq\dots\supseteq R(\AAA/\SSS)
$$
is a decreasing sequence of graded ideals.
Moreover, if $\AAA'$ with the automorphism $\SSS'$ is another category and there is an equivalence $F:\AAA\rightarrow\AAA'$ such that $F\SSS\cong\SSS'F$, then $F$ induces an isomorphism of graded algebras $\ph_F:\ZZ^*(\AAA)\rightarrow \ZZ^*(\AAA')$ such that $\ph_F\big(K_r(\AAA/\SSS)\big)=K_r(\AAA'/\SSS')$.
\end{coro}
\begin{Proof} Suppose that $\upsilon:\SSS'F\rightarrow F\SSS$ is a natural isomorphism. It induces a natural isomorphism $\upsilon^n:(\SSS')^nF\rightarrow F\SSS^n$ for any $n\in\mathbb{Z}$. Then we define $F_{\SSS}:\AAA/\SSS\rightarrow \AAA'/\SSS'$ by the equality $F_{\SSS}x=Fx$ on objects and by the equality $F_{\SSS}(f)=(\upsilon^n)_y^{-1}F(f)\in \big((\AAA'/\SSS')_{Fx}^{Fy}\big)_n$ on morphisms $f\in \big((\AAA/\SSS)_x^y\big)_n=\AAA_x^{\SSS^ny}$. It is easy to see that $F_{\SSS}$ is a degree preserving equivalence.
Now the stated result follows from Theorem \ref{ideals} and Lemma \ref{gc}.
\end{Proof}

\section{Calabi-Yau categories}

In this section we recall the definition of a (weakly) Calabi-Yau category and establish some dualities arising for such categories. We will give an alternative definition of K\"ulshammer ideals for a Calabi-Yau category and show that K\"ulshammer ideals in such categories satisfy additional properties.

\begin{Def}{\rm
Let $\AAA$ be a category with a fixed automorphism $\SSS$ and let $d$ be some integer. The category $\AAA$ is called {\it $d$-Calabi-Yau} if $\AAA_x^y$ is finite dimensional for any $x,y\in\AAA$ and there is a family of linear maps $tr_x:\AAA_x^{\SSS^dx}\rightarrow\kk$ ($x\in\AAA$) such that
\begin{itemize}
\item the pairing $(,):\AAA_y^{\SSS^dx}\times \AAA_x^y\rightarrow\kk$ given by the formula $(f,g)=tr_x(fg)$ is nondegenerate, and
\item for all $m\in\mathbb{Z}$, $g\in\AAA_x^{\SSS^my}$, and $f\in\AAA_y^{\SSS^{d-m}x}$ one has
\begin{equation}\label{CY}
tr_x(\SSS^m(f)g)=(-1)^{m(d-m)}tr_y(\SSS^{d-m}(g)f).
\end{equation}
\end{itemize}
If the first condition is fulfilled and the second condition is true for $m=0$, then $\AAA$ is called a {\it weakly $d$-Calabi-Yau} category. Usually Calabi-Yau categories are assumed to be triangulated and $\SSS$ is assumed to be the shift functor. For more details on Calabi-Yau categories see \cite{Kel}.
}
\end{Def}

If $\AAA$ is a (weakly) $d$-Calabi-Yau category, then we define the map $tr_{\AAA}:\oplus_{x\in\AAA,n\in\bZ}\AAA_x^{\SSS^nx}\rightarrow\kk$ by the equality
$$tr_{\AAA}|_{\AAA_x^{\SSS^nx}}=\begin{cases}
tr_x,&\mbox{ if $n=d$,}\\
0&\mbox{otherwise.}
\end{cases}$$

\begin{Lemma}
Let $\AAA$ with the automorphism $\SSS$ be a weakly $d$-Clabi-Yau category. Then the map $tr_{\AAA}$ induces a map $tr_{\AAA}:\Nat_*(\AAA)\rightarrow\kk$. If $\AAA$ is $d$-Calabi-Yau, then $tr_{\AAA}$
induces also a map $tr_{\AAA}:\Ab_*(\AAA)\rightarrow\kk$.
\end{Lemma}
\begin{Proof} If $\AAA$ is weakly Calabi-Yau, then $tr_{\AAA}(fg)=tr_{\AAA}(\SSS^n(g)f)$ for  $f\in\AAA_x^{\SSS^{n}y}$ and $g\in\AAA_y^{x}$ by definition. Thus, the first assertion is true. For the second assertion it suffices to show that $tr_{\AAA}(f)=(-1)^{n+1}tr_{\AAA}(\SSS(f))$ for all $f\in\AAA_x^{\SSS^nx}$. If $n\not=d$, then both sides are zero. For $n=d$, let us take $m=1$, $y=\SSS^{-1}x$, $g=\Id_x$. By \eqref{CY}, we have $$tr_{\AAA}(\SSS(f))=tr_x(\SSS(f))=(-1)^{d-1}tr_y(f)=(-1)^{d+1}tr_{\AAA}(f)$$ for all $f\in\AAA_{\SSS^{-1}x}^{\SSS^{d-1}x}$. Replacing $x$ by $\SSS x$, we deduce the required equality.
\end{Proof}

Let us now recall the definitions of Hochschild-Mitchel homology and cohomology. Since in this paper we use this notion only for a graded bimodule over a nongraded category, we restrict our definition only to this case.

\begin{Def}{\rm Let $\AAA$ be a nongraded category (i.e. a graded category with all morphisms of degree zero) and let $\MM$ be a graded $\AAA$-bimodule. For $n\ge 0$, we define the set of {\it $n$-cochains}
$$
C^n(\AAA,\MM)=\prod\limits_{x_0,\dots, x_n\in\AAA}\Hom\big(\AAA^{x_0}_{x_1}\ot\cdots\ot \AAA^{x_{n-1}}_{x_n},\MM(x_0,x_n)\big)
$$
and the set of {\it $n$-chains}
$$
C_n(\AAA,\MM)=\bigoplus\limits_{x_0,\dots, x_n\in\AAA}\MM(x_n,x_0)\ot \AAA^{x_0}_{x_1}\ot\cdots\ot \AAA^{x_{n-1}}_{x_n}.
$$
For $n\ge 0$, let us define the linear maps $d_{C^n(\AAA,\MM)}:C^{n}(\AAA,\MM)\rightarrow C^{n+1}(\AAA,\MM)$ and $d_{C_n(\AAA,\MM)}:C_{n+1}(\AAA,\MM)\rightarrow C_n(\AAA,\MM)$ by the equalities
\begin{multline*}
d_{C^n(\AAA,\MM)}(\alpha)(f_0\ot\cdots\ot f_n)=f_0\alpha(f_1\ot\cdots\ot f_n)+\sum\limits_{i=1}^n(-1)^i\alpha(f_0\ot\cdots\ot f_{i-1}f_i\ot \cdots\ot f_n)\\+(-1)^{n+1}\alpha(f_0\ot \cdots\ot f_{n-1})f_n,\\
d_{C_n(\AAA,\MM)}(u\ot f_0\cdots\ot f_{n})=uf_0\ot f_1\ot\cdots\ot f_{n}+\sum\limits_{i=1}^n(-1)^iu\ot f_0\ot\cdots\ot f_{i-1}f_i\ot \cdots\ot f_n\\+(-1)^{n+1}f_nu\ot f_0\ot \cdots\ot f_{n-1},
\end{multline*}
where $\alpha\in C^n(\AAA,\MM)$, $u\in\MM(x,y)$ for some $x,y\in\AAA$ and $f_0,\dots,f_n$ are morphisms in $\AAA$ such that the composition $f_0\dots f_n$ is a well-defined element of $\AAA^y_x$. We define the $n$-th {\it Hochschild-Mitchel cohomology and homology} of $\AAA$ with coefficients in $\MM$ by the equalities
$$\HH^n(\AAA,\MM)=\Ker d_{C^n(\AAA,\MM)}/\Im d_{C^{n-1}(\AAA,\MM)}\mbox{ and }\HH_n(\AAA,\MM)=\Ker d_{C_{n-1}(\AAA,\MM)}/\Im d_{C_n(\AAA,\MM)}.$$
Here we set for convenience $d_{C^{-1}(\AAA,\MM)}=d_{C_{-1}(\AAA,\MM)}=0$.
}
\end{Def}

\begin{rema} For more details on Hochschild-Mitchel cohomology and homology see \cite{Mit}, derived invariance of Hochschild-Mitchell homology and cohomology is proved in \cite{HS}. In general, the formulas for the differentials above have to be more complicated (see, for example, \cite{TamTs}), but in the case where $\AAA$ is nongraded the formulas above are valid. In fact, in this case we can even define $\HH_n(\AAA,\MM)$ and $\HH^n(\AAA,\MM)$ for a nongraded $\MM$ and then set 
$\HH_n(\AAA,\MM)=\oplus_{i\in\mathbb{Z}}\HH_n(\AAA,\MM_i)$ and $\HH^n(\AAA,\MM)=\oplus_{i\in\mathbb{Z}}\HH^n(\AAA,\MM_i)$ for a graded bimodule $\MM$. 
\end{rema}

As usually, we have $\HH^0(\AAA,\MM)=\MM^{\AAA}$ and $\HH_0(\AAA,\MM)=\MM_{\AAA}$. Thus, by Lemmas \ref{gc} and \ref{gab} we have $\Nat^*(\AAA,\AAA/\SSS)\cong\HH^0(\AAA/\SSS)$ and $\Nat_*(\AAA)\cong\HH_0(\AAA,\AAA/\SSS)$ for a category $\AAA$ with an automorphism $\SSS$.

If $\MM$ is a graded $\AAA$-linear category, then the sets defined above carry a lot of additional structure. Here we need only the so-called contraction map $$\ii:C^m(\AAA,\MM)\rightarrow \Hom\big(C_n(\AAA,\MM), C_{n-m}(\AAA,\MM)\big).$$
In general $\ii$ is defined for all integers $n\ge m$ and induces a well-defined map from $\HH^m(\AAA,\MM)$ to $\Hom\big(\HH_n(\AAA,\MM), \HH_{n-m}(\AAA,\MM)\big)$ that is denoted by $\ii$ as well. In the present paper we need only the case $n=m$. Let us write $\ii_{\alpha}$ for the image of $\alpha\in C^n(\AAA,\MM)$ in $\Hom\big(C_n(\AAA,\MM), C_0(\AAA,\MM)\big)$ under the map $\ii$. Then $\ii:C^n(\AAA,\MM)\rightarrow \Hom\big(C_n(\AAA,\MM), C_0(\AAA,\MM)\big)$ can be defined by the equality
$$
\ii_{\alpha}(u\ot f_1\ot\cdots\ot f_n)=u\circ \alpha(f_1\ot\cdots\ot f_n)
$$
for $\alpha\in C^n(\AAA,\MM)$, $u\in\MM(x,y)$ and morphisms $f_1,\dots,f_n$ in $\AAA$ such that the composition $f_1\dots f_n$ is a well-defined element of $\AAA^y_x$, where $x,y\in\AAA$. As was mentioned above, $\ii$ induces a well defined map
$\ii:\HH^n(\AAA,\MM)\rightarrow \Hom\big(\HH_n(\AAA,\MM), \HH_0(\AAA,\MM)\big)$, i.e. we have a map $\ii_{\alpha}:\HH_n(\AAA,\MM)\rightarrow \MM_{\AAA}$ for each $\alpha\in\HH^n(\AAA,\MM)$.

Suppose now that $\SSS$ is an automorphism of $\AAA$ that acts on the graded bimodule $\MM$. Then there is an actions of $\SSS$ on $C^n(\AAA,\MM)$ and $C_n(\AAA,\MM)$ defined by the equalities
$$
({}^{\SSS}\alpha)(f_1\ot\cdots\ot f_n)={}^{\SSS}\big(\alpha({}^{\SSS^{-1}}\!\!f_1\ot\cdots\ot {}^{\SSS^{-1}}\!\!f_n)\big)\mbox{ and }{}^{\SSS}(u\ot f_1\ot\cdots\ot f_n)={}^{\SSS}\!u\ot {}^{\SSS}\!f_1\ot\cdots\ot {}^{\SSS}\!f_n.
$$
It is easy to see that this action induces an action of $\SSS$ on $\HH^n(\AAA,\MM)$ and $\HH_n(\AAA,\MM)$, respectively. As before, $a\SSS$ acts on $C_n(\AAA,\MM)$ and $\HH_n(\AAA,\MM)$ for any $a\in\kk^*$. Moreover, $\ii$ induces maps
$$
\ii:C^n(\AAA,\MM)^{\SSS}\rightarrow \Hom\big(C_n(\AAA,\MM)_{a\SSS}, C_0(\AAA,\MM)_{a\SSS}\big)
$$
and
$$
\ii:\HH^n(\AAA,\MM)^{\SSS}\rightarrow \Hom\big(\HH_n(\AAA,\MM)_{a\SSS}, (\MM_{\AAA})_{a\SSS}\big).
$$

There are gradings on $C^n(\AAA,\MM)$ and $C_n(\AAA,\MM)$ defined by the equalities $C^n(\AAA,\MM)_i=C^n(\AAA,\MM_{i})$ and $C_n(\AAA,\MM)_i=C_n(\AAA,\MM_{i})$. These gradings induce gradings on $\HH^n(\AAA,\MM)$, $C^n(\AAA,\MM)^{\SSS}$, $\HH^n(\AAA,\MM)^{\SSS}$, $\HH_n(\AAA,\MM)$, $C_n(\AAA,\MM)_{a\SSS}$, and $\HH_n(\AAA,\MM)_{a\SSS}$. If $\VV$ is a graded space, then we define $\VV^*$ as a graded space with the degree $i$ component $(\VV^*)_i=(\VV_{-i})^*$ for $i\in\mathbb{Z}$. If we fix some isomorphism  $\Theta:\XX\rightarrow \YY^*$ between graded spaces $\XX$ and $\YY$, then, for $\UU\subseteq \YY$, we define $\UU^{\perp}=\{x\in \XX\mid \Theta(x)|_{\UU}=0\}$. If $\UU\subseteq \YY$ is a graded subspace, then it is easy to see that $\Theta$ induces an isomorphism from $\UU^{\perp}$ to  $(\YY/\UU)^*$.

Since $(\AAA/\SSS)_{\AAA}\cong\Nat_*(\AAA)$ and $\Big((\AAA/\SSS)_{\AAA}\Big)_{-\SSS}\cong \Ab_*(\AAA)$, we have a map $tr_{\AAA}:(\AAA/\SSS)_{\AAA}\rightarrow\kk$ in the case where $\AAA$ is weakly Calabi-Yau and a map $tr_{\AAA}:\Big((\AAA/\SSS)_{\AAA}\Big)_{-\SSS}\rightarrow\kk$ in the case where $\AAA$ is Calabi-Yau.

\begin{theorem}\label{dual} Let $\AAA$ with the automorphism $\SSS$ be a weakly $d$-Calabi-Yau category. Then, for any $n\ge 0$, the map $\Theta_n:\HH^n(\AAA,\AAA/\SSS)\rightarrow \HH_n(\AAA,\AAA/\SSS)^*[-d]$ defined by the equality
$\Theta_n(\alpha)=tr_{\AAA}\ii_{\alpha}$
for $\alpha\in\HH^n(\AAA,\AAA/\SSS)$ is an isomorphism of graded spaces. Moreover, if $\AAA$ is $d$-Calabi-Yau, then $\Theta_n$ induces an isomorphism of graded spaces $\Theta_n:\HH^n(\AAA,\AAA/\SSS)^{\SSS}\rightarrow \big(\HH_n(\AAA,\AAA/\SSS)_{-\SSS}\big)^*[-d]$.
\end{theorem}
\begin{Proof} Let us consider the map $\Theta_n:C^n(\AAA,\AAA/\SSS)\rightarrow C_n(\AAA,\AAA/\SSS)^*[-d]$ defined by the equality $\Theta_n(\alpha)=tr_{\AAA}\ii_{\alpha}$. Its $i$-th component equals to the composition of isomorphisms
\begin{multline*}
\prod\limits_{x_0,\dots, x_n\in\AAA}\Hom\big(\AAA^{x_0}_{x_1}\ot\cdots\ot \AAA^{x_{n-1}}_{x_n},\AAA^{\SSS^{i}x_0}_{x_n}\big) \cong\prod\limits_{x_0,\dots, x_n\in\AAA}\Big(\AAA^{\SSS^{i}x_0}_{x_n}\ot \Hom(\AAA^{x_0}_{x_1}\ot\cdots\ot \AAA^{x_{n-1}}_{x_n},\kk)\Big)\\
\cong\prod\limits_{x_0,\dots, x_n\in\AAA}\Big((\AAA^{\SSS^{d-i}x_n}_{x_0})^*\ot (\AAA^{x_0}_{x_1}\ot\cdots\ot \AAA^{x_{n-1}}_{x_n})^*\Big)\cong \left(\bigoplus\limits_{x_0,\dots, x_n\in\AAA}\AAA^{\SSS^{d-i}x_n}_{x_0}\ot \AAA^{x_0}_{x_1}\ot\cdots\ot \AAA^{x_{n-1}}_{x_n}\right)^*,
\end{multline*}
where the isomorphism $\AAA^{\SSS^{i}x_0}_{x_n}\cong(\AAA^{\SSS^{d-i}x_n}_{x_0})^*$ is induced by the pairing $(,):\AAA^{\SSS^{d-i}x_n}_{x_0}\times\AAA^{\SSS^{i}x_0}_{x_n} \rightarrow\kk$ defined by the equality $(f,g)=tr_{x_n}\big(\SSS^{i}(f)g\big)$.
Thus, $\Theta_n$ is an isomorphism.

Let us now prove that $\Theta_{n+1}d_{C^n(\AAA,\AAA/\SSS)}=(d_{C_n(\AAA,\AAA/\SSS)})^*\Theta_{n}$ for any $n\ge 0$. Fix some $\alpha\in C^n(\AAA,\MM)_i$, $x,y\in\AAA$, $u\in\AAA/\SSS(x,y)_j$, and some morphisms $f_0,\dots,f_n$  in $\AAA$ such that the composition $f_0\dots f_n$ is a well-defined element of $\AAA^y_x$.
Then $$\big(\Theta_{n+1}d_{C^n(\AAA,\AAA/\SSS)}(\alpha)\big)(u\ot f_0\ot\cdots\ot f_n)=tr_{\AAA}\big(ud_{C^n(\AAA,\AAA/\SSS)}(\alpha)(f_0\ot\cdots\ot f_n)\big)$$
and
$$\big((d_{C_n(\AAA,\AAA/\SSS)})^*\Theta_{n}(\alpha)\big)(u\ot f_0\ot\cdots\ot f_n)=tr_{\AAA}\ii_{\alpha}d_{C_n(\AAA,\AAA/\SSS)}(u\ot f_0\ot\cdots\ot f_n).$$
It is easy to see that
\begin{multline*}
ud_{C^n(\AAA,\AAA/\SSS)}(\alpha)(f_0\ot\cdots\ot f_n)-\ii_{\alpha}d_{C_n(\AAA,\AAA/\SSS)}(u\ot f_0\ot\cdots\ot f_n)\\=(-1)^{n+1}(u\alpha(f_0\ot\cdots\ot f_{n-1})f_n-f_nu\alpha(f_0\ot\cdots\ot f_{n-1})).
\end{multline*}
If $i+j\not=d$, then it follows from the definition of $tr_{\AAA}$ that $$tr_{\AAA}\big(u\alpha(f_0\ot\cdots\ot f_{n-1})f_n\big)=tr_{\AAA}\big(f_nu\alpha(f_0\ot\cdots\ot f_{n-1})\big)=0.$$ For $i+j=d$, we have $$tr_{\AAA}\big(u\alpha(f_0\ot\cdots\ot f_{n-1})f_n\big)=tr_{\AAA}\big(f_nu\alpha(f_0\ot\cdots\ot f_{n-1})\big)$$ by the weak Calabi-Yau property. Thus, the maps $\Theta_i$ induce an isomorphism of graded spaces $$\Theta_n:\HH^n(\AAA,\AAA/\SSS)\cong \Ker (d_{C_n(\AAA,\MM)})^*[-d]/\Im (d_{C_{n-1}(\AAA,\MM)})^*[-d]\cong\HH_n(\AAA,\AAA/\SSS)^*[-d]$$ for each $n\ge 0$. So the first part of the theorem is proved.

Suppose now that $\AAA$ is $d$-Calabi-Yau. Note that it follows from the definition of $tr_{\AAA}$ that $\Theta_n(\alpha)(\gamma)=0$ for $\alpha\in C^n(\AAA,\AAA/\SSS)_i$ and $\gamma\in C_n(\AAA,\AAA/\SSS)_j$ if $i+j\not=d$. Then, for $\alpha\in C^n(\AAA,\AAA/\SSS)$, $x,y\in\AAA$, $u\in\AAA/\SSS(x,y)$, and morphisms $f_1,\dots,f_n$  in $\AAA$ such that the composition $f_1\dots f_n$ is a well-defined element of $\AAA^y_x$, we have
\begin{multline*}
\Theta_n(\alpha)\big(u\ot f_1\ot\dots f_n+{}^{\SSS}(u\ot f_1\ot\dots f_n)\big)=tr_{\AAA}\big(u\alpha(f_1\ot\dots \ot f_n)+{}^{\SSS}u\alpha({}^{\SSS}f_1\ot\dots {}^{\SSS}f_n)\big)\\
=tr_{\AAA}\big(u(\alpha-{}^{\SSS^{-1}}\!\!\alpha)(f_1\ot\dots \ot f_n)\big)=\Theta_n(\alpha-{}^{\SSS^{-1}}\!\!\alpha)(u\ot f_1\ot\dots \ot f_n).
\end{multline*}
This holds by the Calabi-Yau property and since $\SSS^{-1}$ acts on $(\AAA/\SSS)_d$ as $(-1)^d\SSS$ acts on $\AAA$.
Let us consider $\UU_n=\{\gamma+{}^{\SSS}\gamma\mid \gamma\in \HH_n(\AAA,\AAA/\SSS)\}\subseteq \HH_n(\AAA,\AAA/\SSS)$. It follows from the argument above that $\Theta_n(\alpha)(\gamma+{}^{\SSS}\gamma)=\Theta_n(\alpha-{}^{\SSS^{-1}}\!\!\alpha)(\gamma)$ for $\alpha\in\HH^n(\AAA,\AAA/\SSS)$ and $\gamma\in \HH_n(\AAA,\AAA/\SSS)$. Since $\Theta_n:\HH^n(\AAA,\AAA/\SSS)\rightarrow \HH_n(\AAA,\AAA/\SSS)^*[-d]$ is an isomorphism, we have $(\UU_n[d])^{\perp}=\HH^n(\AAA,\AAA/\SSS)^{\SSS}$, where $\perp$ is defined with respect to $\Theta_n$. Then, by the discussion above, $\Theta_n$ induces an isomorphism $$\Theta_n:\HH^n(\AAA,\AAA/\SSS)^{\SSS}\cong (\HH_n(\AAA,\AAA/\SSS)/\UU_n)^*[-d]=\big(\HH_n(\AAA,\AAA/\SSS)_{-\SSS}\big)^*[-d].$$
\end{Proof}

In particular, it follows from Theorem \ref{dual} that if $\AAA$ is weakly $d$-Calabi-Yau, then there is an isomorphism $\Theta_{tr}=\Theta_0:\Nat^*(\AAA)\rightarrow \big(\Nat_*(\AAA)\big)^*[-d]$ and, if $\AAA$ is $d$-Calabi-Yau, then $\Theta_{tr}$ induces an isomorphism $\Theta_{tr}:\ZZ^*(\AAA)\rightarrow \big(\Ab_*(\AAA)\big)^*[-d]$. Whenever we consider some Calabi-Yau category $\AAA$, we fix the isomorphism $\Theta_{tr}$ and write $(f,g)$ instead of $\Theta_{tr}(f)(g)$ for $f\in \ZZ^*(\AAA)$ and $g\in \Ab_*(\AAA)$.

By Lemma \ref{gab}, the map $\xi_p$ from Lemma \ref{ppower} induces a map $\xi_p:\Ab_*(\AAA)\rightarrow \Ab_*(\AAA)$. Recall that $T_r(\AAA/\SSS)$ is the kernel of the map $\xi_p^r$.

\begin{Lemma}\label{CYdef}
Let $\AAA$ with the automorphism $\SSS$ be a $d$-Calabi-Yau category. Then $K_r(\AAA/\SSS)=K_{r,d}(\AAA/\SSS)=T_r(\AAA/\SSS)^{\perp}$.
\end{Lemma}
\begin{Proof} It is clear that $K_r(\AAA/\SSS)\subseteq K_{r,d}(\AAA/\SSS)\subseteq T_r(\AAA/\SSS)^{\perp}$. Let us prove that $T_r(\AAA/\SSS)^{\perp}\subseteq K_r(\AAA/\SSS)$. Suppose that $\eta\in \ZZ^n(\AAA)$ does not belong to $K_r(\AAA/\SSS)$. Then there is $f\in \Ab_m(\AAA)$ such that $f^{p^r}=0$ and $\eta f\not=0$. Then there exists $\theta\in \ZZ^{d-m-n}(\AAA)$ such that 
$$0\not=(\theta,\eta f)=tr_{\AAA}(\theta\eta f)=(-1)^{n(d-m-n)}tr_{\AAA}(\eta\theta f)=(-1)^{n(d-m-n)}(\eta,\theta f).$$
Since $\theta f\in T_r(\AAA/\SSS)$, we have $\eta\not\in T_r(\AAA/\SSS)^{\perp}$. Consequently, $K_r(\AAA/\SSS)=T_r(\AAA/\SSS)^{\perp}$.
\end{Proof}

Lemma \ref{CYdef} gives an alternative definition of K\"ulshammer ideals for a Calabi-Yau category. In particular, it implies that $T_r(\AAA/\SSS)^{\perp}$ does not depend on the choice of $\Theta_{tr}$. As usually, such a definition can be reformulated in terms of the adjoint map.

\begin{prop}
Let $\AAA$ with the automorphism $\SSS$ be a $d$-Calabi-Yau category. Assume that the field $\kk$ is perfect. Then, for each $r\ge 0$, there is a unique linear map $\zeta_r:\ZZ^*(\AAA)\rightarrow \ZZ^*(\AAA)$ such that $(\zeta_rf,g)^{p^r}=(f,\xi_p^rg)$ for all $f\in \ZZ^*(\AAA)$ and $g\in \Ab_*(\AAA)$. Moreover, $K_r(\AAA/\SSS)=\Im\zeta_r$.
\end{prop}
\begin{Proof} Let  $\phi_r:\kk\rightarrow \kk$ denote the Frobenius  automorphism $\phi_r(a)=a^{p^r}$.  Then it is easy to see that the map $(\phi_r^{-1})_*(\xi_p^r)^*:\big(\Ab_*(\AAA)\big)^*\rightarrow \big(\Ab_*(\AAA)\big)^*$ defined by the equality $(\phi_r^{-1})_*(\xi_p^r)^*(u)=\phi_r^{-1}u\xi_p^r$ for $u\in\big(\Ab_*(\AAA)\big)^*$ is linear. The equality that $\zeta_r$ has to satisfy can be rewritten in the form $\Theta_{tr}\zeta_r=(\phi_r^{-1})_*(\xi_p^r)^*\Theta_{tr}$.
Thus, the unique map satisfying the required condition is $\zeta_r=\Theta_{tr}^{-1}(\phi_r^{-1})_*(\xi_p^r)^*\Theta_{tr}$. Note that $(\phi_r^{-1})_*(\xi_p^r)^*(u)=0$ for $u\in(\Ab_n(\AAA))^*$ if $p^r\nmid n$ and hence $\zeta_r(f)=0$ for $f\in \ZZ^m(\AAA)$ if $p^r\nmid d-m$.

Let us now prove that $K_r(\AAA/\SSS)=\Im\zeta_r$. It is clear that $\Im\zeta_r\subseteq (\Ker\xi_p^r)^{\perp}$. Let us prove the inverse inclusion. Note that $\Theta_{tr}(\Im\zeta_r)=\Im\big((\phi_r^{-1})_*(\xi_p^r)^*\Theta_{tr}\big)=\Im((\phi_r^{-1})_*(\xi_p^r)^*\big)$. Let us prove that 
$\Theta_{tr}\big((\Ker\xi_p^r)^{\perp}\big)\subseteq\Im((\phi_r^{-1})_*(\xi_p^r)^*\big)$. Suppose that $$h\in\Theta_{tr}\big((\Ker\xi_p^r)^{\perp}\big)_m\subseteq \big(\Ab_{d-m}(\AAA)\big)^*.$$
Let $u_i$, $i\in I\cup J$ be a basis of $\Ab_{d-m}(\AAA)$ such that $u_i$, $i\in J$ is a basis of $\Ker \xi_p^r|_{\Ab_{d-m}(\AAA)}$. Then $\xi_p^r(u_i)$, $i\in I$ is a set of linearly independent elements and there exists $u\in \big(\Ab_{(d-m)p^r}(\AAA)\big)^*$ such that $u\xi_p^r(u_i)=\phi_rh(u_i)$. Then
$h=(\phi_r^{-1})_*(\xi_p^r)^*u\in\Im((\phi_r^{-1})_*(\xi_p^r)^*\big)$. Hence, $K_r(\AAA/\SSS)=\Im\zeta_r$ and the proposition is proved.
\end{Proof}

\section{K\"ulshammer ideals in the Hochschild cohomology}

In this section we will define K\"ulshammer ideals in the Hochschild cohomology of an algebra. From here on $\Lambda$ denotes a $\kk$-algebra.

Let $\DD \Lambda$ be the derived category of the module category over $\Lambda$ and $\Kb \Lambda$ be the full subcategory of $\DD \Lambda$ formed by objects  isomorphic to complexes of finitely generated projective modules.
Let $[1]:\DD \Lambda\rightarrow \DD \Lambda$ be the shift functor. Then we can define the graded center of $\Kb \Lambda$ as the graded center $\Kb \Lambda$ with the automorphism $[1]$.
Let us recall the definition of the homomorphism $\chi_{\Lambda}:\HH^*(\Lambda)\rightarrow \ZZ^*(\Kb \Lambda)$, which is called the {\it characteristic homomorphism}.
It is well known that there is an isomorphism of algebras $\HH^*(\Lambda)\cong \oplus_{n\ge 0}\Hom_{\DD(\Lambda^{op}\ot \Lambda)}(\Lambda,\Lambda[n])$ and so any element of $\HH^n(\Lambda)$ corresponds to a unique morphism $f\in\Hom_{\DD(\Lambda^{op}\ot \Lambda)}(\Lambda,\Lambda[n])$.
Then, for each $X\in \Kb \Lambda$, we define $$\chi_{\Lambda}(f)_X=\Id_X\otimes_{\Lambda}^Lf:X\cong X\otimes_{\Lambda}^L\Lambda\rightarrow X\otimes_{\Lambda}^L\Lambda[n]\cong X[n].$$
It is not hard to see that $\chi_{\Lambda}(f)$ is a natural transformation satisfying the equality $\chi_{\Lambda}(f)_{X[1]}=(-1)^n\chi_{\Lambda}(f)_X[1]$. Thus, $\chi_{\Lambda}:\HH^*(\Lambda)\rightarrow \ZZ^*(\Kb \Lambda)$ is a homomorphism of graded algebras.

It is well known that $\HH^*(\Lambda)$ is invariant under derived equivalences. Since $\Kb \Lambda$ is invariant under derived equivalences too, it is clear that $\ZZ^*(\Kb \Lambda)$ is a derived invariant. In fact, we can say a little more.
Let us recall that by \cite{Ricard} if $\Lambda$ and $\Gamma$ are derived equivalent algebras over a field, then there exist $U\in \DD(\Lambda^{op}\ot \Gamma)$ and $V\in\DD(\Gamma^{op}\ot \Lambda)$ such that $U\ot_{\Gamma}^LV\cong \Lambda$ in $\DD(\Lambda^{op}\ot \Lambda)$ and $V\ot_{\Lambda}^LU\cong \Gamma$ in $\DD(\Gamma^{op}\ot \Gamma)$. In this case $F_U=-\ot_{\Lambda}^LU$ and $F_V=-\ot_{\Gamma}^LV$ induce a pair of quasi inverse equivalences between $\Kb \Lambda$ and $\Kb \Gamma$. Moreover, $U$ and $V$ induce an isomorphism $\ph_{U,V}:\HH^*(\Lambda)\rightarrow \HH^*(\Gamma)$ in the following way. For $f\in\Hom_{\DD (\Lambda^{op}\ot\Lambda)}(\Lambda,\Lambda[n])\cong\HH^*(\Lambda)$ we define $\ph_{U,V}(f)=\alpha_{U,V}^{-1}[n]\big(\Id_V\ot_{\Lambda}^L f\ot_{\Lambda}^L \Id_U\big)\alpha_{U,V}$, where $\alpha_{U,V}:\Gamma\rightarrow V\ot_{\Lambda}^LU$ is an isomorphism and we use the identifications $V\ot_{\Lambda}^L U\cong V\ot_{\Lambda}^L\Lambda\ot_{\Lambda}^L U$ and $V\ot_{\Lambda}^L U[n]\cong V\ot_{\Lambda}^L\Lambda[n]\ot_{\Lambda}^L U$. Note also that due to Corollary \ref{ideals_aut} the equivalence $F_U$ induces an isomorphism $\ph_{F_U}:\ZZ^*(\Kb \Lambda)\rightarrow \ZZ^*(\Kb \Gamma)$.

\begin{Lemma}\label{inv}
Let $\Lambda$, $\Gamma$, $U$ and $V$ be as above. Then $\ph_{F_U}\chi_{\Lambda}=\chi_{\Gamma}\ph_{U,V}$.
\end{Lemma}
\begin{Proof} For $X\in\Kb \Gamma$, we define $$\alpha_X=\Id_X\ot_{\Gamma}^L\alpha_{U,V}:X\cong X\ot_{\Gamma}^L\Gamma\rightarrow X\ot_{\Gamma}^L V\ot_{\Lambda}^L U=F_UF_V(X).$$
Then $\alpha:\Id_{\Kb \Gamma}\rightarrow F_UF_V$ is a natural isomorphism. Using the construction of $\ph_{F_U}$ from Theorem \ref{ideals}, we get, for $f\in\Hom_{\DD(\Lambda^{op}\ot\Lambda)}(\Lambda,\Lambda[n])\cong\HH^*(\Lambda)$ and $X\in\Kb \Gamma$, that
\begin{multline*}
\big(\ph_{F_U}\chi_{\Lambda}(f)\big)_X=\big(\Id_X\ot_{\Gamma}^L\alpha_{U,V}^{-1}[n]\big)\big(\chi_{\Lambda}(f)_{X\ot_{\Gamma}^LV}\ot_{\Lambda}^L \Id_U\big)\big(\Id_X\ot_{\Gamma}^L\alpha_{U,V}\big)\\
=\big(\Id_X\ot_{\Gamma}^L\alpha_{U,V}^{-1}[n]\big)\big(\Id_X\ot_{\Gamma}^L\Id_V\ot_{\Lambda}^Lf\ot_{\Lambda}^L \Id_U\big)\big(\Id_X\ot_{\Gamma}^L\alpha_{U,V}\big)\\
=\Id_X\ot_{\Gamma}^L\Big(\alpha_{U,V}^{-1}[n]\big(\Id_V\ot_{\Lambda}^Lf\ot_{\Lambda}^L \Id_U\big)\alpha_{U,V}\Big)=\big(\chi_{\Gamma}\ph_{U,V}(f)\big)_X.
\end{multline*}
\end{Proof}

Lemma \ref{inv} immediately implies derived invariance of the following ideal in the Hochschild cohomology.

\begin{coro}
The ideal $\Ker(\chi_{\Lambda})\subseteq\HH^*(\Lambda)$ is invariant under derived equivalences.
\end{coro}

For any $s\in\mathbb{Z}$, Corollary \ref{ideals_aut} gives a decreasing sequence of ideals
$$
\ZZ^*(\Kb \Lambda)= K_{0,s}(\Kb \Lambda/[1])\supseteq K_{1,s}(\Kb \Lambda/[1])\supseteq\dots\supseteq K_{r,s}(\Kb \Lambda/[1])\supseteq\dots\supseteq R_s(\Kb \Lambda/[1]).
$$

\begin{Def}{\rm
The $(r,s)$-th {\it higher K\"ulshammer} ideal $\HK_{r,s}^*(\Lambda)$, the $r$-th {\it  higher K\"ulshammer} ideal $\HK_r^*(\Lambda)$, the $s$-th {\it higher Reynolds ideal} $\HR_s^*(\Lambda)$ and the {\it higher Reynolds ideal} $\HR^*(\Lambda)$ of $\Lambda$ are the ideals in $\HH^*(\Lambda)$ defined by the equalities
\begin{multline*}
\HK_{r,s}^*(\Lambda)=\chi_{\Lambda}^{-1}\big(K_{r,s}(\Kb \Lambda/[1])\big),\,\,\HK_r^*(\Lambda)=\chi_{\Lambda}^{-1}\big(K_r(\Kb \Lambda/[1])\big),\\
\HR_s^*(\Lambda)=\chi_{\Lambda}^{-1}\big(R_s(\Kb \Lambda/[1])\big),\mbox{ and }\HR^*(\Lambda)=\chi_{\Lambda}^{-1}\big(R(\Kb \Lambda/[1])\big).
\end{multline*}
}
\end{Def}

\begin{theorem}\label{KulInHH}
If $\Lambda$ is an algebra and $s$ is an integer, then
$$
\HH^*(\Lambda)= \HK_{0,s}^*(\Lambda)\supseteq \HK_{1,s}^*(\Lambda)\supseteq\dots\supseteq \HK_{r,s}^*(\Lambda)\supseteq\dots\supseteq \HR_s^*(\Lambda)
$$
is a decreasing sequence of graded ideals.
In particular,
$$
\HH^*(\Lambda)= \HK_0^*(\Lambda)\supseteq \HK_1^*(\Lambda)\supseteq\dots\supseteq \HK_r^*(\Lambda)\supseteq\dots\supseteq \HR^*(\Lambda)
$$
is a decreasing sequence of graded ideals. Moreover, if $\Gamma$ is derived equivalent to $\Lambda$, then there is an isomorphism of graded algebras $\ph:\HH^*(\Lambda)\rightarrow\HH^*(\Gamma)$ such that $\ph\big(\HK_{r,s}^*(\Lambda)\big)=\HK_{r,s}^*(\Gamma)$.
\end{theorem}
\begin{Proof} $\HK_{r,s}^*(\Lambda)$ is an ideal, since it is an inverse image of the ideal $K_{r,s}(\Kb \Lambda/[1])$ under the algebra homomorphism $\chi_{\Lambda}$. All inclusions follow from the corresponding inclusions for the ideals in $\ZZ^*(\Kb \Lambda)$.

Let now $U\in \DD(\Lambda^{op}\ot \Gamma)$ and $V\in\DD(\Gamma^{op}\ot \Lambda)$ be as above. By Lemma \ref{inv} and Corollary \ref{ideals_aut}, we have $$\chi_{\Gamma}\ph_{U,V}\big(\HK_{r,s}^*(\Lambda)\big)=\ph_{F_U}\chi_{\Lambda}\big(\HK_{r,s}^*(\Lambda)\big)\subseteq\ph_{F_U}\big(K_{r,s}(\Kb \Lambda/[1])\big)=K_{r,s}(\Kb \Gamma/[1]),$$
i.e. $\ph_{U,V}\big(\HK_{r,s}^*(\Lambda)\big)\subseteq\HK_{r,s}^*(\Gamma)$, and
$$\chi_{\Lambda}\ph_{U,V}^{-1}\big(\HK_{r,s}^*(\Gamma)\big)=\ph_{F_U}^{-1}\chi_{\Gamma}\big(\HK_{r,s}^*(\Gamma)\big)\subseteq\ph_{F_U}^{-1}\big(K_{r,s}(\Kb \Gamma/[1])\big)=K_{r,s}(\Kb \Lambda/[1]),$$
i.e. $\ph_{U,V}^{-1}\big(\HK_{r,s}^*(\Gamma)\big)\subseteq\HK_{r,s}^*(\Lambda)$. Thus, $\ph_{U,V}$ satisfies the required conditions as desired.
\end{Proof}

\section{Zero degree}

In this section we fix some algebra $\Lambda$ over a field $\kk$ of characteristic $p$. We will be considering the finite dimensional case, but the first construction is valid for any algebra. We define $[\Lambda,\Lambda]=\{ab-ba\mid a,b\in \Lambda\}$.
It is well known and can be proved analogously to Lemma \ref{ppower} that the map $\xi_p:\Lambda\rightarrow \Lambda$ defined by the equality $\xi_p(a)=a^{p}$ for $a\in \Lambda$ induces a well defined map $\xi_p:\Lambda/[\Lambda,\Lambda]\rightarrow \Lambda/[\Lambda,\Lambda]$, i.e. an endomorphism of $\HH_0(\Lambda)$.

Let us recall the notion of the so-called {\it Hattori-Stallings trace} (see \cite{Hat, Stal, Len}). For the map $f\in\End_\Lambda(\Lambda^n)$, we define $tr(f)=\sum\limits_{i=1}^n\pi_if(e_i)$, where $e_i\in \Lambda^n$ is an element that has $1$ in the $i$-th component and zeros in all others, and $\pi_i:\Lambda^n\rightarrow \Lambda/[\Lambda,\Lambda]$ is the composition of the canonical projection to the $i$-th component and the canonical projection $\Lambda\twoheadrightarrow \Lambda/[\Lambda,\Lambda]$. For a finitely generated projective module $P$ and $f\in\End_{\Lambda}(P)$, we choose some pair of maps $\iota:P\hookrightarrow \Lambda^n$ and $\pi:\Lambda^n\twoheadrightarrow P$ such that  $\pi\iota=1_P$ and define $tr(f)=tr(\iota f\pi)$. One can check that this definition does not depend on $\iota$ and $\pi$. Now, for a bounded complex $C$ with finitely generated projective terms and a map $f\in\End_{\Kb \Lambda}(C)$, we define $tr(f)=\sum_{i\in\mathbb{Z}}(-1)^itr(f_i)$, where $f_i$ is the $i$-th component of $f$. Among other $tr(f)$ has the following properties: $tr(f+h)=tr(f)+tr(h)$, $tr(fh)=tr(hf)$. One can check that $tr:\oplus_{x\in\Kb \Lambda}(\Kb \Lambda)_x^x\rightarrow \Lambda/[\Lambda,\Lambda]$ is a well-defined map that, moreover, induces a map $tr:\Ab_0(\Kb \Lambda)\rightarrow \Lambda/[\Lambda,\Lambda]$. One can easily check also that there is a well defined map $\phi:\Lambda/[\Lambda,\Lambda]\rightarrow \Ab_0(\Kb \Lambda)$ that sends the class of $a\in \Lambda$ to the element $\phi(a)\in \oplus_{x\in\Kb \Lambda}(\Kb \Lambda)_x^x$ that has only one nonzero component $\phi(a)_{\Lambda}:\Lambda\rightarrow \Lambda$ defined by the equality $\phi(a)_{\Lambda}(1_{\Lambda})=a$. Note that $\Lambda/[\Lambda,\Lambda]$ and $\Ab_0(\Kb \Lambda)$ are $\ZZ(\Lambda)$-bimodules, where the second $\ZZ(\Lambda)$-bimodule structure is induced by the inclusion $\chi_{\Lambda}|_{\ZZ(\Lambda)}:\ZZ(\Lambda)\hookrightarrow \ZZ^0(\Kb \Lambda)$.

\begin{Lemma} The maps $\phi$ and $tr$ are homomorphisms of $\ZZ(\Lambda)$-bimodules such that $\xi_p\phi=\phi\xi_p$, $tr\xi_p=\xi_p tr$ and $tr\phi=\Id_{\HH_0(\Lambda)}$.
\end{Lemma}
\begin{Proof} All the assertions can be easily verified.
\end{Proof}

Let us now describe $\Ab_0(\Kb \Lambda)$ for a finite dimensional algebra $\Lambda$. This will allow us to obtain alternative descriptions of $\HK_{r,0}^0(\Lambda)$ and $\HR^0_0(\Lambda)$.
Let us consider the set of indecomposable objects $U\in\Kb \Lambda$ that are not isomorphic to direct summands of $\Lambda$ and satisfy the condition $\max\limits_{H_i(U)\not=0} i=0$.
Let us choose one object in each isomorphism class contained in this set and denote by $C$ the obtained collection of objects of $\Kb \Lambda$.

\begin{theorem}\label{ab0}
Let $\Lambda$ be a finite dimensional algebra.
Then $\Ab_0(\Kb \Lambda)=\Im\phi\oplus\bigoplus\limits_{x\in C}\kk\overline{1}_x$, where $\overline{1}_x$ denotes the class of $\Id_x$ in $\Ab_0(\Kb \Lambda)$. In particular, $\Ker tr=\bigoplus\limits_{x\in C}\kk\big(\overline{1}_x-\phi tr(\overline{1}_x)\big)$.
\end{theorem}
\begin{Proof} Let us prove that $\Im\phi\oplus\bigoplus\limits_{x\in C}\kk\overline{1}_x$ is really a subspace of $\Ab_0(\Kb \Lambda)$, i.e. that if $f+\sum\limits_{i=1}^ka_i\overline{1}_{x_i}=0$ for some $f\in\Im\phi$, $a_i\in\kk$ and distinct $x_i\in C$ ($1\le i\le k$), then $f=0$ and $a_i=0$ for all $1\le i\le k$.
Since $\phi$ is injective, it is enough to show that $\sum\limits_{i=1}^ka_i\overline{1}_{x_i}\not\in\Im\phi$ if at least one of the elements $a_i$ is nonzero. We may assume that $a_1\not=0$. Let us consider $I\subseteq \oplus_{x\in \Kb \Lambda}(\Kb \Lambda)_x^x$ generated by all the nilpotent maps and all the maps that can be factored throw an indecomposable element not isomorphic to an element of the form $x_1[i]$ ($i\in\mathbb{Z}$). Then it is easy to see that any element $g$, whose class in $\Ab_0(\Kb \Lambda)$ coincides with the class of $\sum\limits_{i=1}^ka_i\overline{1}_{x_i}$, can be represented in the form
$$
g=\sum_{x\cong x_1,i\in\mathbb{Z}}c_{x,i}\Id_{x[i]}+T
$$
for some $T\in I$ and $c_{x,i}\in\kk$ almost all zero such that $\sum_{x\cong x_1,i\in\mathbb{Z}}(-1)^ic_{x,i}=a_0$. It is easy to see that $g\not=0$.

Now, note that any element of $\Ab_0$ can be represented by $(f_x)_{x\in \Kb \Lambda}\in \oplus_{x\in \Kb \Lambda}(\Kb \Lambda)_x^x$ such that $f_x=0$ for $x\not\in C\cup\{\Lambda\}$. Indeed, if $x=y\oplus z$, then, for any $f:x\rightarrow x$, we have
$$\overline{f}=\overline{f(\pi_y\iota_y+\pi_z\iota_z)}=\overline{\iota_yf\pi_y}+\overline{\iota_zf\pi_z},$$
where $\pi_y:x\twoheadrightarrow y$, $\pi_z:x\twoheadrightarrow z$, $\iota_y:y\hookrightarrow x$, $\iota_z:z\hookrightarrow x$ are the canonical projections and inclusions, and $\overline{a}$ denotes the class of $a$ in $\Ab_0(\Kb \Lambda)$. Thus, we may assume that $f_x=0$ for any decomposable $x$. Then, due to the equality $\overline{f}+\overline{f[1]}=0$, we may assume that $f_x=0$ if $\max\limits_{H_i(x)\not=0} i\not=0$. Finally, for any $x$ such that $f_x$ is still nonzero, we can choose an isomorphism $\alpha:x\cong y$ ($y\in C$) or a direct inclusion $\alpha:x\hookrightarrow \Lambda$ and change $f_x$ by $\alpha f_x\beta$, where $\beta$ is such a map that $\beta\alpha=\Id_x$.

The class of any element $f\in(\Kb \Lambda)_{\Lambda}^{\Lambda}$ is obviously contained in $\Im\phi$. Let us now take $U\in C$. We may assume that the differential $d_U$ of $U$ has image contained in $UJ_{\Lambda}$, where $J_{\Lambda}$ is the Jacobson radical of $\Lambda$.
Since $\End_{\Kb \Lambda}(U)$ is a local algebra, any $f\in (\Kb \Lambda)_U^U$ can be represented in the form $f=a_U\Id_U+f_N$, where $f_N$ is nilpotent, $a_U \in \kk$. Thus, it remains to show that $\overline{f}\in\Im\phi$ for any nilpotent $f\in (\Kb \Lambda)_U^U$.
Since we assume that $\Im d_U\subseteq UJ_\Lambda$, it is easy to show that all the components of $f$ are nilpotent. Let us prove that $\overline{f}\in\Im\phi$ using induction on the length of $U$. The assertion is obvious if $U$ has only one nonzero term.
Suppose that the assertion holds for complexes of length $n$ and $U$ has length $n+1$, i.e.
$$
U=(\cdots\rightarrow 0\rightarrow U_{-n}\xrightarrow{d_{1-n}} U_{1-n}\xrightarrow{d_{2-n}}\cdots\xrightarrow{d_{-1}}U_{-1}\xrightarrow{d_0}U_0\rightarrow 0\rightarrow\cdots),
$$
and $f$ has components $f_{-n},\dots,f_0$. Let us prove by induction that the class of $f$ in $\Ab_0(\Kb \Lambda)$ equals to the class of the endomorphism of $U(i)$ with components $f_{-n},\dots,f_0$, where
$$
U(i)=(\cdots\rightarrow 0\rightarrow U_{-n}\xrightarrow{d_{1-n}f^i_{-n}} U_{1-n}\xrightarrow{d_{2-n}}\cdots\xrightarrow{d_{-1}}U_{-1}\xrightarrow{d_0}U_0\rightarrow 0\rightarrow\cdots).
$$
The assertion is  vacuous for $i=0$. For the induction step, it is enough to represent the map with components $f_{-n},\dots,f_0$ from $U(i-1)$ to itself as the composition of the map with components $\Id_{U_{-n}},f_{1-n},\dots,f_0$ from $U(i-1)$ to $U(i)$ and the map with components
$f_{-n},\Id_{U_{1-n}},\dots,\Id_{U_0}$ from $U(i)$ to $U(i-1)$. Since the map $f_{-n}$ is nilpotent,
$$U(i)=U_{-n}[n]\oplus (\cdots\rightarrow 0\rightarrow U_{1-n}\xrightarrow{d_{2-n}}U_{2-n}\xrightarrow{d_{3-n}}\cdots\xrightarrow{d_{-1}}U_{-1}\xrightarrow{d_0}U_0\rightarrow 0\rightarrow\cdots)$$
for big enough $i$, and the induction hypothesis implies $\overline f\in\Im\phi$.
\end{Proof}

\begin{Def}{\rm
The $r$-th {\it K\"ulshammer ideal} $K_r\Lambda$ ($n\ge 0$) is the set of such $a\in \ZZ(\Lambda)$ that $ab\in [\Lambda,\Lambda]$ for all $b\in \Lambda$ such that $b^{p^r}\in[\Lambda,\Lambda]$. The {\it Reynolds ideal} of $\Lambda$ is the set $R\Lambda=\cap_{r\ge 0}K_r\Lambda$.
}
\end{Def}

It is easy to see that
$$
\ZZ(\Lambda)= K_0\Lambda\supseteq K_1\Lambda\supseteq\dots\supseteq K_r\Lambda\supseteq\dots\supseteq R\Lambda
$$
is a decreasing sequence of ideals. It is also not difficult to prove that $$R\Lambda=\{a\in \ZZ(\Lambda)\mid aJ_\Lambda\subseteq [\Lambda,\Lambda]\}$$ if $\Lambda$ is finite dimensional, where $J_{\Lambda}$ is the Jacobson radical of $\Lambda$.
Now we are ready to describe the ideals $\HK_{r,0}^0\subseteq \HH^0(\Lambda)=\ZZ(\Lambda)$.

\begin{coro}\label{KIinv}
If $\Lambda$ is a finite dimensional algebra, then $\HK_{r,0}^0(\Lambda)=K_r\Lambda$ for any $r\ge 0$. In particular, $\HR_0^0(\Lambda)=R\Lambda$ in this case.
\end{coro}
\begin{Proof} Since $\phi$ is injective and respects $\xi_p$, we have by definition $$K_r\Lambda=\Ann_{\ZZ(\Lambda)}\Ker(\xi_p^{\Lambda})^r=\Ann_{\ZZ(\Lambda)}(\Ker(\xi_p^{\Ab})^r\cap \Im\phi) \text{ and }$$ $$ \HK_{r,0}^0(\Lambda)=\Ann_{\ZZ(\Lambda)}\Ker(\xi_p^{\Ab})^r$$ where $\xi_p^{\Lambda}$ denotes $\xi_p:\Lambda/[\Lambda,\Lambda]\rightarrow \Lambda/[\Lambda,\Lambda]$ and $\xi_p^{\Ab}$ denotes $\xi_p:\Ab_0(\Kb \Lambda)\rightarrow \Ab_0(\Kb \Lambda)$. By Theorem \ref{ab0} and since $\xi_p^{\Ab}(\Im\phi)\subseteq\Im\phi$ we have $\Ker(\xi_p^{\Ab})^r\subseteq \Im\phi$, and hence the assertion follows.
\end{Proof}

\begin{coro}\label{invGenKul}
If $\Lambda$ and $\Gamma$ are derived equivalent finite dimensional algebras, then there exists an isomorphism $\ph:\ZZ(\Lambda)\cong \ZZ(\Gamma)$ such that $\ph(K_r\Lambda)=K_r\Gamma$ for any $r\ge 0$ and $\ph(R\Lambda)=R\Gamma$.
\end{coro}
\begin{Proof} Follows from Theorem \ref{KulInHH} and Corollary \ref{KIinv}.
\end{Proof}

\begin{rema}
It follows from Corollary \ref{invGenKul} that the set $R\Lambda=\{a\in \ZZ(\Lambda)\mid aJ_{\Lambda}\subseteq [\Lambda,\Lambda]\}$ is an ideal in $\ZZ(\Lambda)$ invariant under derived equivalences if $\Lambda$ is an algebra over a field of characteristic $p$. In fact, our argument can be adopted to prove the derived invariance of this ideal for a finite dimensional algebra over a field of characteristic $0$ as well. For this one uses the fact, that $R\Lambda$ is the annihilator of the set of nilpotent elements in $\Ab_0(\Kb\Lambda)$.
\end{rema}

Let us recall the classical definition of K\"ulshammer ideals.

\begin{Def}{\rm
The finite dimensional algebra $\Lambda$ is  called {\it symmetric} if there is a nondegenerate bilinear form $(,):\Lambda\times \Lambda\rightarrow\kk$ such that $(ab,c)=(a,bc)$ and $(a,b)=(b,a)$ for all $a,b,c\in \Lambda$.
}
\end{Def}

For a symmetric algebra $\Lambda$, the $r$-th {\it classical K\"ulshammer ideal} $K^{cl}_r\Lambda$ ($r\ge 0$) is the set of  $a\in \Lambda$ such that $(a,b)=0$ for all $b\in \Lambda$ such that $b^{p^r}\in[\Lambda,\Lambda]$. It is known that
$$
\ZZ(\Lambda)= K^{cl}_0\Lambda\supseteq K^{cl}_1\Lambda\supseteq\dots\supseteq K^{cl}_r\Lambda\supseteq\dots\supseteq R\Lambda
$$
is a decreasing sequence of ideals. Moreover, if the algebra $\Gamma$ is derived equivalent to $\Lambda$, then $\Gamma$ is symmetric \cite{Ricard} and there is an isomorphism from $\ZZ(\Lambda)$ to $\ZZ(\Gamma)$ that maps $K_r^{cl}\Lambda$ to $K_r^{cl}\Gamma$ for any $r\ge 0$ \cite{Zim1}. The later fact can be recovered from the following lemma, the proof is analogous to the proof of Lemma \ref{CYdef}:

\begin{Lemma}\label{class}
If $\Lambda$ is symmetric, then $K^{cl}_r\Lambda=K_r\Lambda$ for any $r\ge 0$.
\end{Lemma}

The following lemma is well known \cite{Ha}.

\begin{Lemma}\label{Calab}
If $\Lambda$ is a symmetric algebra, then $\Kb \Lambda$ is a $0$-Calabi-Yau category.
\end{Lemma}

\begin{coro}
If $\Lambda$ is symmetric, then $\HK_r^0(\Lambda)=K_r\Lambda$ for any $r\ge 0$. In particular, $\HR^0(\Lambda)=R\Lambda$ is the socle of the algebra $\Lambda$. 
\end{coro}
\begin{Proof} Follows from Lemmas \ref{CYdef} and \ref{class}, and Corollary \ref{KIinv}.
\end{Proof}

\begin{rema}
Let $\Lambda$ be symmetric. The bilinear form $(,):\Lambda\times \Lambda\rightarrow\kk$ induces a nondegenerate bilinear form $(,):\HH^*(\Lambda)\times\HH_*(\Lambda)$ by the equality $(f,u)=\ee(f\frown u)$. Here the map $\ee:\Lambda/[\Lambda,\Lambda]\cong\HH_0(\Lambda)\rightarrow\kk$ is induced by the map $\ee:\Lambda\rightarrow\kk$ defined by the equality
$\ee(a)=(1,a)$. Thus, we can define the map $\lambda_{\Lambda}:\Ab_*(\Kb \Lambda)\rightarrow\HH_*(\Lambda)$ as the unique map satisfying the equality $(f,\lambda_{\Lambda}(u))=(\chi_{\Lambda}(f),u)$ for all $f\in\HH^*(\Lambda)$ and $u\in \Ab_*(\Kb \Lambda)$. It is not difficult to show that $\lambda_{\Lambda}$ is a homomorphism of graded $\HH^*(\Lambda)$-modules. Note that $\Ker(\chi_{\Lambda})\subseteq\HR^*(\Lambda)$. Actually we have $\Ker(\chi_{\Lambda})=\Im(\lambda_{\Lambda})^{\perp}$ and $\HK_r^*(\Lambda)=\lambda_{\Lambda}(T_r\Kb \Lambda)^{\perp}$. 
\end{rema}

\begin{rema} One can consider the Tate-Hochschild cohomology $\widehat\HH^*(\Lambda)$ of a selfinjective algebra $\Lambda$ and the stable category $\umod \Lambda$ of the category of finitely generated $\Lambda$-modules. Note that $\umod \Lambda$ is a triangulated category with the shift functor $\Omega_{\Lambda}^{-1}$.
In this case there exists the characteristic map $\underline\chi_{\Lambda}:\widehat\HH^*(\Lambda)\rightarrow \ZZ^*(\umod \Lambda)$ and one can define K\"ulshamer and Reynolds ideals in $\widehat\HH^*(\Lambda)$ as preimages of the corresponding ideals in $\ZZ^*(\umod \Lambda)$.
These ideals are invariant under stable equivalences of Morita type and it would be interesting to study their properties. Of course, the case where $\Lambda$ is stably $d$-Calabi-Yau for some integer $d$ is of special interest. Note, in particular, that symmetric algebras are stably $(-1)$-Calabi-Yau.
\end{rema}

\section{Example}

In this section we are going to compute all notions defined above for the category $\Kb \Lambda$, where $\Lambda=k[x]/x^2$. The classification of indecomposable objects in $\Kb \Lambda$ is given in \cite{Ku}, the graded center of $\Kb \Lambda$ is computed in \cite{KY}.

For any $m\leq n \in \mathbb{Z}$ consider $$\Lambda^{[m,n]} = \cdots \rightarrow 0 \rightarrow \Lambda \xrightarrow{x} \Lambda \rightarrow \cdots \rightarrow \Lambda \xrightarrow{x} \Lambda \rightarrow 0 \cdots,$$ where nonzero entries are concentrated in the interval $[m,n]$.  Each indecomposable object of $\Kb \Lambda$ is isomorphic to an object of the form $\Lambda^{[m,n]}$ for some $m\leq n \in \mathbb{Z}$. Let us denote by $x_t$ the element of $\Hom(\Lambda^{[m,n]},\Lambda^{[m',n']})$ given by the multiplication by $x$ in degree $t$ and zero in all other degrees. It is easy to see that $\Hom(\Lambda^{[m,n]},\Lambda^{[m,n]})$ is two dimensional and any map is homotopic to a map of the form $c_1\Id_{\Lambda^{[m,n]}}+c_2x_n$ for some $c_1,c_2\in\kk$. If $[m,n]\neq[m',n']$ and the intervals have a nontrivial intersection, then $\Hom(\Lambda^{[m,n]},\Lambda^{[m',n']})$ is one dimensional in the following cases:\\
1) $m\leq m', n\leq n'$, any map is homotopic to a map of the form $cx_n$ ($c\in\kk$):
$$
\xymatrix {
\cdots \ar[r]&0 \ar[r]  & \Lambda \ar[r]^x &\Lambda \ar[r] \ar[d]^0 & \cdots \ar[r] & \Lambda \ar[r]^x &\Lambda \ar[r]^x \ar[d]^0 &\Lambda \ar[r]\ar[d]^{cx}& 0 \ar[r] & \cdots & \\
 &  \cdots \ar[r]&0 \ar[r]  & \Lambda \ar[r]^x &\Lambda \ar[r] & \cdots \ar[r] & \Lambda \ar[r]^x &\Lambda \ar[r]^x& \Lambda \ar[r] & 0 \ar[r] & \cdots   \\
}
$$
2) $m\geq m', n\geq n'$, any map is homotopic to a map of the form $c\Id_{[m,n']}$ ($c\in\kk$), where the map $\Id_{[m,n']}$ is induced by $\Id_{\Lambda}$ in degrees $m, \cdots,n'$ and zero in other degrees:
$$
\xymatrix {
 &  \cdots \ar[r]&0 \ar[r]  & \Lambda \ar[r]^x \ar[d]^{c} &\Lambda \ar[r] & \cdots \ar[r] & \Lambda \ar[r]^x \ar[d]^{c} &\Lambda \ar[r]^x \ar[d]^{c}& \Lambda \ar[r] & 0 \ar[r] & \cdots   \\
\cdots \ar[r]&0 \ar[r]  & \Lambda \ar[r]^x &\Lambda \ar[r]  & \cdots \ar[r] & \Lambda \ar[r]^x &\Lambda \ar[r]^x  &\Lambda \ar[r]& 0 \ar[r] & \cdots & \\
}
$$

In all other cases there are no nonzero morphisms.

The following description of $\Nat^t(\Kb \Lambda)$ and $\ZZ^t(\Kb \Lambda)$ for $t\geq 0$ was obtained in \cite{KY}. 

$\Nat^0(\Kb \Lambda)$ consists of natural transformations $\eta$ given by the data of the form $\{\mu, \lambda_{[m,n]}\in\kk,-\infty<m\le n<\infty\}$, the corresponding natural transformation is given by $\eta_{\Lambda^{[m,n]}}=\mu \Id_{\Lambda^{[m,n]}}+\lambda_{[m,n]} x_n$. 

The transformation $\eta$ belongs to $\ZZ^0(\Kb \Lambda)$ if and only if $\lambda_{[m,n]}=\lambda_{[m+r,n+r]}$ for any $r$.

$\Nat^t(\Kb \Lambda),$ $t>0$ consists of natural transformations $\eta$ given by the data of the form $\{c\in{\bf k}\}$, the corresponding natural transformation is given by $\eta_{\Lambda^{[m,n]}}=c \Id_{[m,n-t]}$ for $n-t-m \geq 0$ and $\eta_{\Lambda^{[m,n]}}=0$, otherwise. 

The transformation $\eta$ belongs to $\ZZ^t(\Kb \Lambda)$ if and only if $\charr\kk=2$ or $t$ is even. 

$\Nat^t(\Kb \Lambda)=0=\ZZ^t(\Kb \Lambda)$ for $t<0$.

$\Nat_0(\Kb \Lambda)=\oplus_{\Lambda^{[m,n]}}\langle \Id_{\Lambda^{[m,n]}} \rangle \oplus V_0$, where $V_0$ is one dimensional. $V_0=\oplus_{\Lambda^{[m,n]}\in\Kb \Lambda}\langle x_n \rangle/U_0$, where $U_0$ is the subspace generated by the elements $x_n - (-1)^{n'-n}x_{n'}$ for $x_n:\Lambda^{[m,n]}\rightarrow \Lambda^{[m,n]}$, $x_{n'}:\Lambda^{[m',n']}\rightarrow \Lambda^{[m',n']}$. 

$\Ab_0(\Kb \Lambda)=\oplus_{\Lambda^{[m,n]}}\langle \Id_{\Lambda^{[m,n]}}\rangle /W_0  \oplus V_0$, where $W_0$ is the subspace generated by the elements $\Id_{\Lambda^{[m,n]}} - (-1)^r\Id_{\Lambda^{[m+r,n+r]}}$. 

$\Nat_t(\Kb \Lambda)=0=\Ab_t(\Kb \Lambda)$, for $t> 0$.

For $t<0$ the space $\Nat_t(\Kb \Lambda) = V_t$ is one dimensional. $V_t=\oplus_{\Lambda^{[m,n]}, n\geq m-t}\langle x_n \rangle/U_t$, where $U_t$ is the subspace generated by the elements $x_n - (-1)^{n'-n}x_{n'}$ for $x_n:\Lambda^{[m,n]}\rightarrow \Lambda^{[m,n]}[t]$, $x_{n'}:\Lambda^{[m',n']}\rightarrow \Lambda^{[m',n']}[t]$. Here, for $t$ odd, the map $x_{l}$ still denotes the multiplication by $x$ in degree $l$ and zero maps in all other degrees.

For $t<0$ the space $\Ab_t(\Kb \Lambda) = V_t$ if $\charr\kk=2$ or $t$ is even, $\Ab_t(\Kb \Lambda) = 0,$ otherwise.

Let us compute the ideals $K_{r,s}(\Kb \Lambda/\SSS)$ and $K_{r}(\Kb \Lambda/\SSS)$.

For any $r>1$ we have

$$T_r:=T_r(\Kb \Lambda/\SSS)=V_0\oplus\bigoplus_{t<0} \Ab_t(\Kb \Lambda).$$

Let us denote by $\tilde{Z}^0$ the subset of $\ZZ^0(\Kb \Lambda)$ given by $\{0,\lambda_{[m,n]}\in\kk,-\infty<m\le n<\infty\}$. 

$$\Ann((T_r)_{t})=\begin{cases}\tilde{Z}^0 \oplus\bigoplus\limits_{l\geq -t+1} \ZZ^l(\Kb \Lambda/\SSS),&\mbox{ if $\charr\kk=2$ or $t$ is even};\\
\oplus_{l} \ZZ^l(\Kb \Lambda/\SSS),&\mbox{ otherwise}.\end{cases}$$

Hence, $R(\Kb \Lambda/\SSS)=K_{r}(\Kb \Lambda/\SSS)=\tilde{Z}^0$.

If $\charr\kk\neq 2$, then
$$
    K_{r,s}(\Kb \Lambda/\SSS)_t=
    \begin{cases}
      0, & \text{for}\ t<0, s\in \mathbb{Z}, \\
      \ZZ^t(\Kb \Lambda/\SSS), & \text{for}\ s>0, t\geq 0, \\
      \ZZ^t(\Kb \Lambda/\SSS), & \text{for}\ (s-t) \text{ odd}, s\leq0, t \geq 0, \\
      \tilde{Z}^0, & \text{for}\ (s-t) \text{ even}, s\leq 0, t=0, \\
      0, & \text{for}\ (s-t) \text{ even}, s\leq 0, t> 0.
    \end{cases}
 $$

If $\charr\kk= 2$, then
$$
    K_{r,s}(\Kb \Lambda/\SSS)_t=
    \begin{cases}
      0, & \text{for}\ t<0, s\in \mathbb{Z}, \\
      \ZZ^t(\Kb \Lambda/\SSS), & \text{for}\ s>0, t\geq 0, \\
      \tilde{Z}^0, & \text{for}\  s\leq 0, t=0, \\
      0, & \text{for}\  s\leq 0, t> 0,
    \end{cases}
 $$

$$R_s(\Kb \Lambda/\SSS)=K_{r,s}(\Kb \Lambda/\SSS).$$

Let us now compute the corresponding ideals in the Hochshild cohomology.
 The bimodule resolution of $\Lambda$ is
$$\cdots \rightarrow \Lambda\otimes\Lambda \xrightarrow{x\otimes1+1\otimes x} \Lambda\otimes\Lambda \xrightarrow{x\otimes1-1\otimes x} \Lambda\otimes\Lambda \xrightarrow{x\otimes1+1\otimes x} \Lambda\otimes\Lambda \xrightarrow{x\otimes1-1\otimes x} \Lambda$$

$$
   \HH^l(\Lambda)=
    \begin{cases}
      \Lambda, & \text{for}\ l=0, \\
      \Lambda/2x\Lambda, & \text{for even}\  l>0, \\
      \Ann(2x), & \text{for odd}\ l.
    \end{cases}
$$

$\HH^0(\Lambda)=\Lambda$, $\chi_{\Lambda}(c+dx)$ is the natural transformation $\eta_{\Lambda^{[m,n]}}$ given by the data $\{c, \lambda_{[m,n]}\}$, where $\lambda_{m,n}=0$ for $m-n$ odd and $\lambda_{m,n}=d$ for $m-n$ even. 

For greater $l$ we are going to compute $\chi_{\Lambda}$ in the following way: any element of $\HH^l(\Lambda)$ gives a map $f$ from the bimodule resolution of $\Lambda$ to its shift, so first we compute $\Lambda^{[m,m+n]}\ot^L_{\Lambda}\Lambda$ as the totalization of a bicomplex, then we compute mutually inverse isomorphisms $\iota_{[m,n]}: \Lambda^{[m,m+n]} \rightarrow \Lambda^{[m,m+n]}\ot^L_{\Lambda}\Lambda$ and $\pi_{[m,n]}: \Lambda^{[m,m+n]}\ot^L_{\Lambda}\Lambda \rightarrow \Lambda^{[m,m+n]}$, then $\chi_{\Lambda}(f)=\pi_{[m,n]}[l](\Id_{\Lambda^{[m,m+n]}}\ot f)\iota_{[m,n]}$. Since the shift does not matter for these computations we can assume $m+n=0.$ Let us use the notation $d_-:=x\otimes1-1\otimes x,$ $d_+:=x\otimes1+1\otimes x.$ As a right module $\Lambda^{op}\ot \Lambda$  is isomorphic to $\Lambda\oplus \Lambda$ (the first $\Lambda$ is generated by $1\ot 1$, the second $\Lambda$ is generated by $x\ot 1$), let us denote by $\iota_1$ the map $\Lambda \xrightarrow{(1,0)^t} \Lambda^{op}\ot \Lambda$, by $\pi_1$ the map $\Lambda^{op}\ot \Lambda \xrightarrow{(1,0)} \Lambda$ and by $x_2$ the map $\Lambda^{op}\ot \Lambda \xrightarrow{(0,x)} \Lambda$. 

The complex $C:=\Lambda^{[-n,0]}\ot^L_{\Lambda}\Lambda$ is the totalization of the bicomplex $\tilde{C}$ with $n$ nonzero rows:

$$
\xymatrix {
& &  \vdots \ar[d] &\vdots \ar[d] & \vdots \ar[d] &\vdots \ar[d]   \\
& \cdots \ar[r]&\Lambda\otimes\Lambda \ar[r]^{d_-} \ar[d]^{-x\ot 1} &\Lambda\otimes\Lambda \ar[r]^{d_+} \ar[d]^{x\ot 1}&\Lambda\otimes\Lambda \ar[r]^{d_-}\ar[d]^{-x\ot 1}  & \Lambda\otimes\Lambda \ar[d]^{x\ot 1} \\
\tilde{C}=& \cdots \ar[r]&\Lambda\otimes\Lambda \ar[r]^{d_-} \ar[d]^{-x\ot 1} &\Lambda\otimes\Lambda \ar[r]^{d_+} \ar[d]^{x\ot 1}&\Lambda\otimes\Lambda \ar[r]^{d_-}\ar[d]^{-x\ot 1}  & \Lambda\otimes\Lambda \ar[d]^{x\ot 1} \\
&\cdots \ar[r]&\Lambda\otimes\Lambda \ar[r]^{d_-}  &\Lambda\otimes\Lambda \ar[r]^{d_+} &\Lambda\otimes\Lambda \ar[r]^{d_-} & \Lambda\otimes\Lambda  \\
}
$$

$$C^{-i}=(\Lambda^{op}\ot \Lambda)^{i+1} \text{ for } i\leq n; C^{-i}=(\Lambda^{op}\ot \Lambda)^{n+1}  \text{ for } i> n.$$ 
To obtain $C^{-i}$ one takes the sum of the entries of the bicomplex along the diagonal, the numbering of the summands $\Lambda^{op}\ot \Lambda$ goes from the lower left entry to the upper right entry. 
The differential $d^{-i}$ is a matrix with entries $d^{-i}_{k,k}=d_{(-1)^{i+k-1}}$, $d^{-i}_{k,k+1}=(-1)^{i+k}x\ot 1$, all other entries are zero. The maps $\iota_{[-n,0]}: \Lambda^{[-n,0]} \rightarrow C$ and $\pi_{[-n,0]}: C \rightarrow \Lambda^{[-n,0]}$ can be defined by the equalities
$$
\begin{aligned}
\iota^{[-n,0]}_{-i}&=((-1)^{\ceil*{\frac{i}{2}}}\iota_1,(-1)^{\ceil*{\frac{i-1}{2}}
}\iota_1, \cdots, \iota_1 )^t\text{ }:\text{ }\Lambda \rightarrow (\Lambda^{op}\ot \Lambda)^{i+1},-n\leq-i\leq0;\\
\pi^{[-n,0]}_{-i}&=((-1)^{\ceil*{\frac{i}{2}}
}\pi_1,0, \cdots, 0 )\text{ }:\text{ } (\Lambda^{op}\ot \Lambda)^{i+1} \rightarrow \Lambda,-n<-i\leq0;\\
\pi^{[-n,0]}_{-n}&=((-1)^{\ceil*{\frac{i}{2}}
}\pi_1\pm x_2,\pm x_2, \cdots, \pm  x_2 ) \text{ }: \text{ } (\Lambda^{op}\ot \Lambda)^{n+1}\rightarrow \Lambda,
\end{aligned}
$$
for appropriate signs before the maps $x_2$.

Case 1: $char\kk \neq 2$


For $l$ odd $\HH^l(\Lambda)$ is generated by $x$. We have $\chi_{\Lambda}\HH^l(\Lambda)=0$ since $\ZZ^l(\Kb \Lambda/\SSS)=0$.

For $l$ even  $\HH^l(\Lambda)$ is generated by the class of $1$. The corresponding $f$ has entries $f_{-i}=\Id_{\Lambda}$, $n\geq i\geq l$. The corresponding composition $(\pi_{[-n,0]}[l](\Id_{\Lambda^{[-n,0]}}\ot f)\iota_{[-n,0]})_{-i}=(-1)^{l/2}\Id_{\Lambda}$ for $n\geq i\geq l$. Hence, $$\chi_{\Lambda}\HH^l(\Lambda)=\ZZ^l(\Kb \Lambda/\SSS), \text{ for } l>0.$$
  If $r>0$, $s\leq 0$ is even, then $\HK_{r,s}(\Lambda)=\langle x \rangle \oplus \bigoplus\limits_{l\ge 0}\HH^{2l+1}(\Lambda)$ since
 $$
    \HK^l_{r,s}(\Lambda)=
    \begin{cases}
      \langle x \rangle, & \text{for}\  l=0, \\
      \HH^l(\Lambda), & \text{for}\ l > 0, l \text{ odd},  \\
      0, & \text{for } l> 0,l \text{ even}.
    \end{cases}
$$
  If $r>0$ and either $s\leq 0$ is odd or $s>0$, then $\HK_{r,s}(\Lambda)=\HH^*(\Lambda).$

  $$\HR(\Lambda)=\HK_r(\Lambda)=\langle x \rangle \oplus \bigoplus\limits_{l\ge 0}\HH^{2l+1}(\Lambda); \HR_s(\Lambda)=\HK_{r,s}(\Lambda), \text{for }r>0.$$
  
Case 2: $char\kk = 2$  
  
For $l>0$, $\HH^l(\Lambda)$ is generated by $1$ and $x$, for $x$ the corresponding $f$ has entries $f_{-i}=x\ot 1$, $n\geq i\geq l$. The corresponding composition $\pi_{[-n,0]}[l](\Id_{\Lambda^{[-n,0]}}\ot f)\iota_{[-n,0]}$ is clearly homotopic to zero, i.e. $\chi_{\Lambda}(f)=0$. For $1$ the corresponding $f$ has entries $f_{-i}=\Id_{\Lambda}$, $n\geq i\geq l$. The corresponding composition $(\pi_{[-n,0]}[l](\Id_{\Lambda^{[-n,0]}}\ot f)\iota_{[-n,0]})_{-i}=\Id_{\Lambda}$ for $n\geq i\geq l$. Hence, $$\chi_{\Lambda}\HH^l(\Lambda)=\ZZ^l(\Kb \Lambda/\SSS), \text{ for } l>0.$$

  If $r>0$, $s\leq 0$, then $\HK_{r,s}(\Lambda)=\bigoplus_{t\geq0}\langle x \rangle_t,$ is the ideal generated by $x\in\HH^0(\Lambda)$. Here $\langle x \rangle_t$  denotes the subspace $\langle x \rangle$ of $\HH^t(\Lambda)$.

  If $r>0$, $s> 0$, then $\HK^*_{r,s}(\Lambda)=      \HH^*(\Lambda)$  
  
$$\HK_r(\Lambda)=\bigoplus_{t\geq0}\langle x \rangle_t; \HR_s(\Lambda)=\HK_{r,s}(\Lambda), \text{for } r>0.$$

\end{document}